\documentclass[notitlepage,leqno,10pt]{article}
\usepackage[pdf4]{pstricks}
\usepackage{amsmath,amssymb,amsthm,enumitem,graphicx,relsize,pst-all,pstricks-add,pst-bspline,float}
\usepackage{amscd}

\newcommand
{\pstEllipse
}[5][]{
\psset
{#1}
\parametricplot
{#4}{#5}{#2\space t cos mul
#3\space t sin mul
}}

\renewcommand{\d}{\delta }
\newcommand{\D }{\Delta }

\renewcommand{\l }{\lambda }

\newcommand{\n }{\nabla }

\newcommand{\s }{\sigma }

\renewcommand{\o }{\omega }

\newcommand{\ov}{\overline}
\newcommand{\intbar}{\mathop{\int\makebox(-13.5,0){\rule[4pt]{.7em}{0.3pt}}%
\kern-6pt}\nolimits}

\newcommand{\be}{\begin{equation}}
\newcommand{\ee}{\end{equation}}
\newcommand{\bes}{\begin{equation*}}
\newcommand{\ees}{\end{equation*}}
\newcommand{\ba}{\begin{eqnarray}}
\newcommand{\ea}{\end{eqnarray}}
\newcommand{\bas}{\begin{eqnarray*}}
\newcommand{\eas}{\end{eqnarray*}}
\newenvironment{pf}{\noindent{\sc Proof}.\enspace}{\rule{2mm}{2mm}\medskip}
\newenvironment{pfn}{\noindent{\sc Proof}}{\rule{2mm}{2mm}\medskip}

\newcommand{\R}{\mathbb{R}}

\newcommand{\Z}{\mathbb{Z}}

\newcommand{\N}{\mathbb{N}}
\renewcommand{\o }{\omega }

\author{Martin MAYER$^{a}$,\;\; Cheikh Birahim NDIAYE$^b$}

\date{}

\title{\bf Barycenter technique and the Riemann mapping problem of Escobar}

\begin{document}

\newtheorem{lem}{Lemma}[section]
\newtheorem{pro}[lem]{Proposition}
\newtheorem{thm}[lem]{Theorem}
\newtheorem{rem}[lem]{Remark}
\newtheorem{cor}[lem]{Corollary}
\newtheorem{df}[lem]{Definition}

\maketitle

\begin{center}

{\small

\noindent  $^{a, b}$ Mathematisches Institut der Justus-Liebig-Universit\"at Giessen, \\Arndtstrasse 2, D-35392 Giessen, Germany.

}
\
\
{\small

\noindent

}

\end{center}

\footnotetext[1]{E-mail addresses: martin.g.mayer@math.uni-giessen.de, ndiaye@everest.mathematik.uni-tuebingen.de, cheikh.ndiaye@math.uni-giessen.de,  Mohameden.Ahmedou@math.uni-giessen.de.}

\

\

\begin{center}
{\bf Abstract}

\end{center}
We solve the remaining cases of the Riemann mapping problem of Escobar\cite{es2}. Indeed, performing a suitable scheme of the barycenter technique of Bahri-Coron\cite{bc} via the Chen\cite{chen}'s bubbles, we solve the cases left open after the work of Chen\cite{chen}. Thus, combining our work with the ones of Almaraz\cite{almaraz1}, Chen\cite{chen}, Escobar\cite{es2},\cite{es4}, and Marques\cite{marques},\cite{marques1}, we have that every compact Riemannian manifold with boundary of dimension greater or equal than \;$3$\; carries a conformal scalar flat metric with constant mean curvature.

\begin{center}

\bigskip\bigskip
\noindent{\bf Key Words:} Scalar curvature, Mean curvature, Variational bubbles, Algebraic topological methods.
\bigskip

\centerline{\bf AMS subject classification:  53C21, 35C60, 58J60, 55N10.}

\end{center}
\section{Introduction and statement of the results}
In his attempt to suitably generalize the celebrated Riemann mapping theorem of complex analysis which asserts that any simply connected proper domain of the plane is conformally diffeomorphic to a disk, Escobar\cite{es2} raised the question of whether every \;$n$-dimensional compact Riemannian manifold with boundary and \;$n\geq 3$\; carries a conformal scalar flat Riemannian metric with constant mean curvature. In \cite{es2} and \cite{es4}, Escobar provides a positive answer when \;$n=3$, $n=4$\; or \;when \;$ n=5$\; and the boundary is umbilic, and  when \;$n\geq 6$\; with the boundary being non umbilic or the Riemannian manifold being locally conformally flat and the boundary being umbilic. Later, Marques\cite{marques},\cite{marques1} gives a positive answer to some remaining cases, precisely when \;$n=4$\; or \;$5$\; and the boundary is not umbilic, when \;$n\geq 8$\; and the boundary is umbilic and  not locally conformally flat with respect to the induced Riemannian metric, and when \;$n\geq 9$\; with the boundary being umbilic and the Weyl tensor does not vanish identically on the boundary. In \cite{almaraz1}, Almaraz\cite{almaraz1} gives a positive answer when \;$n=6, 7, 8$, the boundary is umbilic, and the Weyl tensor does not vanish identically on the boundary. Recently, Chen\cite{chen} resolves the problem for many situations of the cases remaining after the above cited works and reduces the other ones to the positivity of the ADM mass of some class of asymptotically flat Riemannian manifolds, like she did in a joint work with S. Brendle for the boundary Yamabe problem in \cite{bs}. However, like in \cite{bs}, the latter positivity is {\em not know to hold}.
\vspace{6pt}

\noindent
Our main goal in this work is to use the algebraic topological argument of Bahri-Coron\cite{bc} to solve the cases left open by Almaraz\cite{almaraz1}, Chen\cite{chen}, Escobar\cite{es2},\cite{es4}, and Marques\cite{marques},\cite{marques1}, like we did in \cite{martndia1} for the boundary Yamabe problem to settle the cases remaining after the works of Escobar\cite{es1} and Brendle-Chen\cite{bs}. Indeed, performing a suitable scheme of the barycenter technique of Bahri-Coron\cite{bc} via the Chen\cite{bs}'s bubbles, we prove a result for the Riemann mapping problem of Escobar which covers all the cases left open after the above cited works. In order to state clearly our theorem, we first fix some notation. Given $(\ov M, g)$ a $n$-dimensional compact Riemannian manifold with boundary \;$\partial M$, interior \;$M$\;and \;$n\geq 3$, we denote by \;$L_g=-4\frac{n-1}{n-2}\D_g+R_g$\; the conformal Laplacian of $(\ov M, g)$ and \;$B_g=\frac{4(n-1)}{n-2}\frac{\partial}{\partial n_g}+2(n-1)H_g$\; the conformal Neumann operator of $(M, g)$, with \;$R_g$\; denoting the scalar curvature of \;$(\ov M, g)$, \;$\D_g$ \;denoting the Laplace-Beltrami operator with respect to $g$, \;$H_g$\; is the mean curvature of \;$\partial M$\; in \;$(\ov M, g)$,  \;$\frac{\partial}{\partial n_g}$ is the outer Neumann operator on $\partial M$ with respect to $g$. Furthermore, we define the following Escobar functional 
\begin{equation}\label{eq:escobarfunctional}
\mathcal{E}_g(u):=\frac{\langle L_gu, u\rangle+\langle B_g u, u\rangle}{(\oint_{\partial M}u^{\frac{2(n-1)}{n-2}}dS_g)^{\frac{n-2}{n-1}}}, \;\;\;\;\;u\in W^{1, 2}_{+}(\ov M),
\end{equation}
where \;$\langle L_gu, u\rangle:=\langle L_gu, u\rangle_{L^2(M)}$, $\langle B_gu, u\rangle:=\langle B_gu, u\rangle_{L^2(\partial M)}$, $dS_g$ is the volume form with respect to the Riemannian metric induced by $g$ on $\partial M$, and $W^{1, 2}_+(\ov M):=\{u\in W^{1, 2}(\ov M): \;\;u>0\}$\; with $W^{1, 2}(\ov M)$\; denoting the usual Sobolev space of functions which are \;$L^2$-integrable with their first derivatives (for more information, see \cite{aubin} and \cite{gt}). Moreover, we recall that the Sobolev quotient of $(M, \partial M, g)$ is defined as 
\begin{equation}\label{eq:Sobolevquotient}
 \mathcal{Q}(M, \partial M, g):=\inf_{u\in W^{1, 2}_+(\ov M)}\mathcal{E}_g(u).
\end{equation}
Now, having fixed the needed notation, we are ready to state our theorem which reads as follows.
\begin{thm}\label{eq:existence}
Assuming that \;$(\ov M, g)$\; is a $n$-dimensional compact Riemannian manifold with boundary \;$\partial M$\; and interior $M$\; such that $\partial M$ is umbilic in $(\ov M, g)$, \;$n\geq 6$, and \;$\mathcal{Q}(M, \partial M, g)>0$,\; then \;$(\ov M, g)$\; carries a conformal scalar flat Riemannian metric with respect to which \;$\partial M$\; has constant mean curvature. 
\end{thm}
\vspace{6pt}

\noindent
Hence, since the only open cases for the Riemann mapping problem of Escobar\cite{es2} is when the dimension of the manifold is greater or equal than \;$6$ with umbilic boundary and positive Sobolev quotient, then clearly Theorem \ref{eq:existence} and the works of Almaraz\cite{almaraz1}, Chen\cite{chen}, Escobar\cite{es2},\cite{es4}, and Marques\cite{marques},\cite{marques1} imply the following positive answer to the high-dimensional generalization by Escobar\cite{es2} of the celebrated Riemann mapping problem of complex analysis.
\begin{thm}\label{eq:BoundaryYamabe}
Every \;$n$-dimensional compact Riemannian manifold with boundary and\;$n\geq 3$\; carries a conformal scalar flat Riemannian metric with respect to which its boundary has constant mean curvature.
\end{thm}
\vspace{10pt}

\noindent
To give a positive answer to the high-dimensional generalization by Escobar\cite{es2} of the celebrated Riemann mapping problem of Riemann surface theory is equivalent to solving a second order elliptic boundary value problem with critical Sobolev nonlinearity on the boundary. Indeed, under the assumptions of Theorem \ref{eq:existence}, the Riemann mapping problem of Escobar\cite{es2} is equivalent to finding a smooth and positive solution of the following semilinear elliptic boundary value problem\\
\begin{equation}\label{eq:bvp}
\left\{
\begin{split}
L_gu&=0 \;\;&\text{in}\;\;M,\\
 B_gu&=2(n-1)u^{\frac{n}{n-2}}\;\;&\text{on}\;\;\partial M.
\
\end{split}
\right.
\end{equation}
The boundary value problem \eqref{eq:bvp} has a variational structure. Indeed, thanks to the work of Cherrier\cite{cherrier}, smooth solutions of \eqref{eq:bvp} can be found by looking at critical points of the Escobar functional \;$\mathcal{E}_g$, and like in \cite{martndia1} we will pursue such an approach here. Precisely, we will perform a suitable application of the barycenter technique of Bahri-Coron\cite{bc} via the Chen\cite{chen}'s bubbles. We describe briefly the barycenter technique of Bahri-Coron\cite{bc} (focusing on \;$\mathcal{E}_g$) for those who are not familiar with. The algebraic topological argument of Bahri-Coron\cite{bc} belongs the the class of indirect methods. Precisely, it is an argument by contradiction. Thus, assuming that the Euler-Lagrange functional \;$\mathcal {E}_g$\; has no critical points, one looks for a contradiction by using the quantization and strong interaction phenomenon that \;$\mathcal{E}_g$\;verifies and the structure of the space of barycenters of \;$\partial M$. To describe how the latter works, we feel more useful for the sake of understanding of reader to do it with figures rather than exact mathematical formulas, of course at the price of precision but with the right intuition of what is going on with the barycenter technique of Bahri-Coron\cite{bc}. First of all, and recalling the assumption \;$\mathcal {E}_g$\; has no critical points (to keep in mind), one has that the quantization phenomenon that \;$\mathcal{E}_g$\; enjoys implies the following figure
\begin{figure}[H]
\begin{center}
\psscalebox{1}
{
\begin{pspicture}(-5.5,0)(9,7)\psset{xunit=25pt,yunit=25pt,runit=20pt}
\put(-8.5,2.2){\large$B_{1}(\partial M)=$}
\psellipse(-2.5,2)(2,0.5)
\put(-4,3.3){$\partial M$}

\psline{->}(0,2)(5,2)
\put(3,3)
{\psscalebox{1 0.8}{\psplot[algebraic,linewidth=0.01,plotpoints=10000]{-1.7}{1.4}{5/(1+25*x^2)}}
}
\put(1.5,0){$+$ Quantization}

\psline(6.5,0)(6.5,8)

\psline(6,1)(7,1)
\put(9.3,1){$W_{0}$}

\psline(6,7)(7,7)
\put(9.3,8.5){$W_{1}$}

\psellipse(9,6)(2,0.5)
\end{pspicture}
}

\end{center}
\caption{}\label{Figure1}
\end{figure}\noindent
which traduces the fact that by bubbling \;$\partial M$\; survives topologically between the first and second critical levels of \;$\mathcal{E}_g$\;that we denote by \;$(W_1, W_0)$ (see \eqref{eq:defenergylevel} for its precise definition). Next, realizing \;$B_2(\partial M)$\;(for its definition see \eqref{eq:barytop}) as a cone over \;$B_1(\partial M)=\partial M$\; with top \;$\partial M$, one has that the quantization phenomenon that \;$\mathcal {E}_g$\; satisfies implies again the following figure
\begin{figure}[H]
\begin{center}
\psscalebox{1}
{
\begin{pspicture}(-4,-3)(9,7)\psset{xunit=25pt,yunit=25pt,runit=20pt}
\put(-8,3){\large$B_{2}(\partial M)=$}
\psline(-5,-2)(-3,6.5)(-1,-2)
\psellipse(-3,-2)(2,0.5)
\put(-4.3,8.5){$\partial M$}
\put(-4.6,-2.7){$B_{1}(\partial M)$}

\psline{->}(-0.5,-2)(6,-2)
\put(3,-2)
{\psscalebox{1 0.7}
{
\psplot[algebraic,linewidth=0.01,plotpoints=10000]{-2.5}{2.5}
{
5/(1+25*(x)^2)
}
}
}
\put(1.5,-4.6){$+$ Quantization}

\psline{->}(-0.5,3)(6,3)
\put(3,4)
{\psscalebox{1 0.7}
{
\psplot[algebraic,linewidth=0.01,plotpoints=10000]{-2.5}{2.5}
{
5/(1+25*(x+1)^2)
+
5/(1+25*(x-1)^2)
}
}
}

\psline(6.5,-3.5)(6.5,8)

\psline(6,-3)(6.75,-3)
\put(8.5,-4){$W_{0}$}

\psline(6,2)(6.75,2)
\put(8.5,2.25){$W_{1}$}

\psline(6,7)(6.75,7)
\put(8.5,8.5){$W_{2}$}

\psline(7.5,1.2)(9.5,5.5)(11.5,1.2)
\psellipse(9.5,1.2)(2,0.5)

\end{pspicture}
}

\end{center}
\caption{}\label{Figure2}
\end{figure}\noindent
which shows the fact that by bubbling \;$B_2(\partial M)$\; as a cone over \;$B_1(\partial M)$\; survives as a nontrivial cone between the second and third critical levels of\; $\mathcal{E}_g$ that we denote by $(W_2, W_1)$. Similarly, realizing\; $B_3(\partial M)$\; as a cone over \;$B_2(\partial M)$\; with top \;$\partial M$, one has that the quantization phenomenon  that \;$\mathcal{E}_g$\;verifies implies again the following figure
\begin{figure}[H]
\begin{center}
\psscalebox{1}
{
\begin{pspicture}(-4,-5)(9,6)\psset{xunit=25pt,yunit=25pt,runit=20pt}
\put(-8,1.5){\large$B_{3}(\partial M)=$}
\psline(-5,-4)(-3,6.5)(-1,-4)
\psline(-5,-4)(-3,0)(-1,-4)
\psline[linestyle=dashed](-3,6.5)(-3,0)
\psellipse(-3,-4)(2,0.5)
\put(-4.3,9){$\partial M$}
\put(-4.7,-3.5){$B_{2}(\partial M)$}

\psline{->}(-0.7,3.5)(6,3.5)
\put(3,4.5)
{\psscalebox{1 0.5}
{
\psplot[algebraic,linewidth=0.01,plotpoints=10000]{-3}{3}
{
5/(1+25*(x)^2)
+
5/(1+25*(x+2)^2)
+
5/(1+25*(x-2)^2)
}
}
}

\psline{->}(-0.7,-0.5)(6,-0.5)
\put(3,-0.5)
{\psscalebox{1 0.5}
{
\psplot[algebraic,linewidth=0.01,plotpoints=10000]{-2}{2}
{
5/(1+25*(x+1)^2)
+
5/(1+25*(x-1)^2)
}
}
}

\psline{->}(-0.7,-4.5)(6,-4.5)
\put(3,-5.5)
{\psscalebox{1 0.5}
{
\psplot[algebraic,linewidth=0.01,plotpoints=10000]{-1}{1}
{
5/(1+25*(x)^2)
}
}
}
\put(1.5,-7){$+$ Quantization}

\psline(6.5,-5.5)(6.5,8)

\psline(6,-5)(6.7,-5)
\put(8.5,-6.5){$W_{0}$}

\psline(6,-1)(6.7,-1)
\put(8.5,-1.5){$W_{1}$}

\psline(6,3)(6.7,3)
\put(8.5,3.5){$W_{2}$}

\psline(6,7)(6.7,7)
\put(8.5,8.5){$W_{3}$}

\psline(7.5,-2)(9.5,4.2)(11.5,-2)
\psline(7.5,-2)(9.5,2)(11.5,-2)
\psline[linestyle=dashed](9.5,2)(9.5,4.2)
\psellipse(9.5,-2)(2,0.5)
\end{pspicture}
}

\end{center}
\caption{}\label{Figure3}
\end{figure}\noindent
which traduces the fact that by bubbling \;$B_3(\partial M)$\; as a cone over \;$B_2(\partial M)$\; survives as a nontrivial cone between the third and fourth critical levels of\; $\mathcal{E}_g$ that we denote by $(W_3, W_2)$. Hence, recursively for $p\in \N^{*}$, realizing\; $B_{p+1}(\partial M)$\; as a cone over \;$B_{p}(\partial M)$\; with top \;$\partial M$, we have that the quantization phenomenon  that \;$\mathcal{E}_g$\;enjoys implies again the following figure
\begin{figure}[H]
\begin{center}
\psscalebox{1}
{
\begin{pspicture}(-4,-3)(9,6)\psset{xunit=25pt,yunit=25pt,runit=20pt}
\put(-8,3){\large$B_{p+1}(\partial M)=$}
\psline(-4.5,-2)(-3,6.5)(-1.5,-2)
\psline[linewidth=0.7pt](-4.5,-2)(-3,1.5)(-1.5,-2)
\psline[linestyle=dashed](-3,6.5)(-3,1.5)
\psellipse(-3,-2)(1.5,0.5)
\put(-4.3,8.5){$\partial M$}
\put(-4.6,-1.2){$B_{p}(\partial M)$}

\psline{->}(-1,-2)(6,-2)
\put(3,-2.25)
{\psscalebox{1 0.6}
{
\psplot[algebraic,linewidth=0.01,plotpoints=10000]{-3.4}{-0.5}
{
5/(1+25*(x+2)^2)
}
\psplot[algebraic,linewidth=0.01,plotpoints=10000]{0.5}{3.4}
{
5/(1+25*(x-2)^2)
}
}
}
\put(2.8,-2.25){$//$}
\put(1.5,-1){$..$ p-times $..$}
\put(1.5,-4.15){$+$ Quantization}

\psline{->}(-1,2.25)(6,2.25)
\put(3,3)
{\psscalebox{1 0.6}
{
\psplot[algebraic,linewidth=0.01,plotpoints=10000]{-3.4}{-0.5}
{
5/(1+25*(x+2)^2)
}
\psplot[algebraic,linewidth=0.01,plotpoints=10000]{0.5}{3.4}
{
5/(1+25*(x-2)^2)
}
}
}
\put(2.8,3){$//$}
\put(1.5,4){$..$ (p+1)-times $..$}

\psline(6.5,-3)(6.5,8)

\psline(6,1)(6.75,1)
\put(8.5,1){$W_{p}$}

\psline(6,7)(6.75,7)
\put(8.5,8.5){$W_{p+1}$}

\psline(7.75,0.6)(9.75,1.3)(11.75,0.6)
\psline[linewidth=0.7pt](7.75,0.6)(9.75,1)(11.75,0.6)
\psline[linewidth=0.7pt,linestyle=dashed,dash=0.5pt](9.75,1.3)(9.75,1)
\rput(9.75,0.6)
{
\pstEllipse[linecolor=black]{2}{0.2}{-180}{0}
}
\rput(9.75,0.6)
{
\pstEllipse[linestyle=dashed,linecolor=black]{2}{0.2}{55}{125}
}

\end{pspicture}
}

\end{center}
\caption{}\label{Figure4}
\end{figure}\noindent
which shows the fact that by bubbling \;$B_{p+1}(\partial M)$\; as a cone over \;$B_{p}(\partial M)$\; survives as a nontrivial cone between the $(p+1)$ and $(p+2)$ critical levels of\; $\mathcal{E}_g$ that we denote by $(W_{p+1}, W_{p})$. On the other hand, the latter recursion leads to a contradiction because of the strong interaction phenomenon  that \;$\mathcal{E}_g$\;enjoys. To see this, we first recall that the \;$\mathcal{E}_g$-energy of the sum of $p$ highly concentrated bubbles $\varphi_{a_i, \l}$ (with $p\in \N$, $p>1$, $a_i\in \partial M$ \;for $i=1, \cdots, p$\; are the concentration points and \;$\l$\; is the concentration parameter) modeling a configuration with optimal weights is roughly given by the sum of self and interaction of the bubbles as follows
\begin{equation}\label{eq:energyexplain}
\mathcal{E}_g(\sum_{i=1}^p \varphi_{a_i, \l})=\sum_{i=1}^p\left(\mathcal{E}_g(\varphi_{a_i, \l})+\sum_{j=1, j\neq i}^pInt_g(\varphi_{a_i, \l}, \varphi_{a_j, \l})\right),
\end{equation}
where \;$Int_g(\varphi_{a_i, \l}, \varphi_{a_j, \l})$\; denotes the interaction of the bubbles $\varphi_{a_i, \l}$\;and $\varphi_{a_j, \l}$. Furthermore, the concentration points live in a world with quantization and strong interaction phenomenon as shown by the following figure
\begin{figure}[H]
\begin{center}
\psscalebox{1}
{
\begin{pspicture}(-4,-0.5)(9,8)\psset{xunit=25pt,yunit=25pt,runit=20pt}
\psellipse[linestyle=dashed](2,5)(5,5)

\def\Arrow
{
\psline{<-}(-2,5)(0,5)
\psline(-3,5)(-2,5)
\psline{>-}(0,5)(2,5)
}

\rput{180}(4,10){\Arrow}
\rput{-60}(-3.3,4.22){\Arrow}
\rput{-120}(-3.8,4.9){\Arrow}

\put(2.5,6.7){$\frac{MD_{g}(a_{i})}{\lambda^{n-2}}$}
\put(-2,12){$\frac{MD_{g}(a_{j_{1}})}{\lambda^{n-2}}$}
\put(-2,0){$\frac{MD_{g}(a_{j_{2}})}{\lambda^{n-2}}$}
\put(9,6){$\frac{MD_{g}(a_{j_{3}})}{\lambda^{n-2}}$}

\put(-1.4,9){$\frac{G_{g}(a_{i},a_{j_{1}})}{\lambda^{n-2}}$}
\put(1.4,3.5){$\frac{G_{g}(a_{i},a_{j_{2}})}{\lambda^{n-2}}$}
\put(5.3,5.5){$\frac{G_{g}(a_{i},a_{j_{3}})}{\lambda^{n-2}}$}


\put(2.52,6.25){\psdot[dotsize=6pt]}
\put(8.75,6.25){\psdot[dotsize=6pt]}
\put(-0.6,0.85){\psdot[dotsize=6pt]}
\put(-0.6,11.65){\psdot[dotsize=6pt]}
\end{pspicture}
}

\end{center}
\caption{}\label{Figure5}
\end{figure}\noindent
with
\begin{equation}\label{eq:selfinteract}
\mathcal{E}_g(\varphi_{a_i, \l})=c_{\infty}-\frac{MD_g(a_i)}{\l^{n-2}}, \;\;\;\text{and}\;\;\;Int_g(\varphi_{a_i, \l}, \varphi_{a_j, \l})=-\frac{G_g(a_i, a_j)}{\l^{n-2}},
\end{equation}
where \;$c_{\infty}$\; is the first critical value, \;$MD_g(a_i)$\; denotes the mass distribution of the particle $\;a_i$\; and \;$G_g$\; is the Green's function of the couple conformal Laplacian and conformal Neumann operator under a suitable normalization. Moreover, \;$MD_g$\; is an \;$L^{\infty}$-function and\; $\min_{\ov M^2}G_g>0$\; which is what we mean by \;$\mathcal{E}_g$\; verifies a strong interaction phenomenon. Thus, clearly \eqref{eq:energyexplain} and \eqref{eq:selfinteract} imply the following figure

\begin{figure}[H]
\begin{center}
\psscalebox{1}
{
\begin{pspicture}(-4,-3)(9,6)\psset{xunit=25pt,yunit=25pt,runit=20pt}
\put(-8,2.5){\large$B_{p_{0}}(\partial M)=$}
\psline(-4.5,-2)(-3,6.5)(-1.5,-2)
\psline(-4.5,-2)(-3,1)(-1.5,-2)
\psline[linestyle=dashed](-3,6.5)(-3,1)
\psellipse(-3,-2)(1.5,0.5)
\put(-4.3,8.5){$\partial M$}
\put(-5,-1.7){$B_{p_{0}-1}(\partial M)$}

\psline{->}(-1,-2)(6,-2)
\put(3,-2.25)
{\psscalebox{1 0.6}
{
\psplot[algebraic,linewidth=0.01,plotpoints=10000]{-3.4}{-0.5}
{
5/(1+25*(x+2)^2)
}
\psplot[algebraic,linewidth=0.01,plotpoints=10000]{0.5}{3.4}
{
5/(1+25*(x-2)^2)
}
}
}
\put(2.8,-2.25){$//$}
\put(1.5,-1){$..$ (p$_{0}$-1)-times $..$}
\put(1.5,-4.15){$+$ Quantization}

\psline{->}(-1,2.25)(6,2.25)
\put(3,3)
{\psscalebox{1 0.6}
{
\psplot[algebraic,linewidth=0.01,plotpoints=10000]{-3.4}{-0.5}
{
5/(1+25*(x+2)^2)
}
\psplot[algebraic,linewidth=0.01,plotpoints=10000]{0.5}{3.4}
{
5/(1+25*(x-2)^2)
}
}
}
\put(2.8,3){$//$}
\put(1.5,4){$..$ p$_{0}$-times $..$}

\psline(6.5,-3)(6.5,8)

\psline(6,1)(6.75,1)
\put(8.5,1){$W_{p_{0}-1}$}

\psline(6,7)(6.75,7)
\put(8.5,8.5){$W_{p_{0}}$}

\put(10,1.5){No space for cone}
\psline[linewidth=3pt](8,1)(11.1,1)

\end{pspicture}
}

\end{center}
\caption{}\label{Figure6}
\end{figure}\noindent
which shows the fact that for \;$p_0$ large, $B_{p_0}(\partial M)$\; as a cone over \;$B_{p_0-1}(\partial M)$\; can not be embedded by bubbling and still be a nontrivial cone between the \;$p_0$\; and \;$p_0+1$\; critical levels of \;$\mathcal{E}_g$. Hence, clearly the figures \eqref{Figure4} and \eqref{Figure6} lead to a contradiction.
\vspace{8pt}

\noindent
The structure of this paper is as follows. In Section \ref{eq:notpre}, we fix some notation and give some preliminaries, like the set of formal barycenters of \;$\partial M$\; and present some useful topological properties of them. Furthermore, we recall the Chen\cite{bs}'s bubbles and the fact that they can be used to replace the standard bubbles in the analysis of diverging Palais-Smale (PS) sequences of the Euler-Lagrange functional \;$\mathcal{E}_g$. Moreover, using a result of Almaraz\cite{almaraz2} and another one of Chen\cite{bc}, we derive self and interaction estimates for the Chen\cite{chen}'s bubbles. In Section \ref{eq:mappingbary}, we use the latter estimates to map the space of barycenter of \;$\partial M$\; of any order into suitable sublevels of \;$\mathcal{E}_g$\; via the Chen\cite{chen}'s bubbles. Finally, in Section \ref{eq:algtop}, we define the neighborhood of potential critical points at infinity of \;$\mathcal{E}_g$ and use the results of Section \ref{eq:mappingbary} to carry our scheme of the barycenter technique of Bahri-Coron\cite{bc} to prove Theorem \ref{eq:existence}.
\vspace{8pt}

\noindent
\begin{center}
{\bf Acknowledgements}
\end{center}
C. B. Ndiaye has been supported  by the DFG project "Fourth-order uniformization type theorems for $4$-dimensional Riemannian manifolds".
\section{Notation and preliminaries}\label{eq:notpre}
In this section, we fix some notations and give some preliminaries. First of all, since the problem under study is conformally invariant and we are dealing with the umbilic case, then from now until the end of the paper \;$(\ov M, g)$\; will be the given underlying compact \;$n$-dimensional Riemannian manifold with boundary \;$\partial M$\; and interior \;$M$, \;$\partial M$\; is totally geodesic in \;$(\ov M, g)$, \;$n\geq 6$, and \;$\mathcal{Q}(M, \partial M, g)>0$. 
\vspace{6pt}

\noindent
In the following, for any Riemannian metric \;$\bar g$\; on \;$\ov M$, we will use the notation \;$B^{\bar g}_{p}(r)$\; to denote the geodesic ball with respect to $\bar g$ of radius \;$r$\;and center \;$p$. Similarly, for $p\in \partial M$, we use the notation \;$\hat {B}^{\bar g}_{p}(r)$\; to denote the geodesic ball in $\partial M$ with respect to the Riemannian metric\;$\hat{\bar g}$\;induced by \;$\bar g$\; on \;$\partial M$\; of radius \;$r$\;and center \;$p$.  We also denote respectively by \;$d_{\bar g}(x,y)$\; the geodesic distance with respect to $\bar g$ between two points \;$x$\;and \;$y$\; of \;$\ov M$ and $d_{\hat{\bar g}}(x, y)$, the geodesic distance with respect to $\hat{\bar g}$ between two points \;$x$\;and \;$y$\; of \;$\partial M$. $inj_{\bar g}(\ov M)$, $inj_{\hat {\bar g}}(\partial M)$\;stand for the injectivity radius of \;$(\ov M, \bar g)$, $( \partial M, \hat{\bar g})$. $dV_{\bar g}$\;denotes the Riemannian measure associated to the metric\;$\bar g$, and $dS_{\bar g}$\; the volume form on \;$\partial M$\; with respect to \;$\hat{\bar g}$\; on \;$\partial M$. For simplicity, we will use respectively \;$B_p(r)$\; and $\hat B_p(r)$ to denote \;$B^g_{p}(r)$\; and \;$\hat B^g_p(r)$. For \;$a\in \ov M$,  we use the notation \;$exp_a^{\bar g}$\; to denote the exponential map with respect to \;$\bar g$\; and set for simplicity $\exp_a:=\exp_a^g$. For \;$a\in \partial M$, we denote by \;$\hat{\exp}_a^{\bar g}$ the exponential map with respect to \;$\hat{\bar g}$, and set \;$\hat{\exp}_a:=\hat{\exp}_a^{g}$.
\vspace{6pt}

\noindent
$\N$\;denotes the set of nonnegative integers, $\N^*$\;stands for the set of positive integers, and  for $k\in \N^*$, $\R^k$\;stands for the standard $k$-dimensional Euclidean space,  $\R^k_+$ the open positive half-space of $\R^k$, and $\bar \R^k_+$ its closure in $\R^k$. For simplicity, we will use the notation \;$\R_+:=\R^1_+$, and $\bar \R_+:=\bar \R^1_+$. For $r>0$, \;$B^{\R^k}_0(r)$ denotes the open ball of \;$\R^k$\; of center \;$0$\; and radius \;$r$, and set for simplicity \;$B^k:=B^{\R^k}_0(1)$. We use \;$g_{B^k}$\; to denote the Euclidean metric on \;$B^k$. For \;$p\in \N^*$, $\sigma_p$\; denotes the permutation group of \;$p$ elements, $(\partial M)^p$ denotes the cartesian product of \;$p$\; copies of \;$\partial M$. For \;$p\in \N^*$, $F((\partial M)^p)$ denotes the fat diagonal of $(\partial M)^p$, namely $F((\partial M)^p):=\{A:=(a_1, \cdots, a_p)\in ( \partial M)^P:\;\;\exists\;i\neq j\;\;\text{with}\;a_i=a_j\}$.  We define $((\partial M)^2)^*:=(\partial M)^2\setminus Diag((\partial M)^2)$ where \;$Diag((\partial M)^2)$\; is the diagonal of $(\partial M)^2$, namely $Diag((\partial M)^2):=\{(a, a): \;\,a\in \partial M\}$. $\sigma_p$ stand for the permutations group of \;$p$\; elements and \;$\D_{p-1}$ the following simplex \;$\D_{p-1}:=\{(\alpha_1, \cdots, \alpha_p):\alpha_i\geq 0, i=1, \cdots, p, \sum_{i=1}^p\alpha_i=1\}$. 
\vspace{6pt}

\noindent
For $1\leq p\leq \infty$ and $k\in \N$, $\beta\in  ]0, 1[$, $L^p(M)$\; and \;$L^p(\partial M)$, \;$W^{k, p}(M)$, $C^k(\ov M)$, and \;$C^{k, \beta}(\ov M)$ stand respectively for the standard $p$-Lebesgue space on \;$M$ and \;$\partial M$, $(k, p)$-Sobolev space, $k$-continuously differentiable space and $k$-continuously differential space of H\"older exponent \;$\beta$, all with respect to \;$g$ (if the definition needs a metric structure) and for precise definitions and properties, see  for example \cite{aubin} or \cite{gt}.
\vspace{6pt}

\noindent
For \;$a\in \partial M$, $O_a(1)$ stands for quantities bounded uniformly in $a$. For\;$\epsilon$\, positive and small, and \;$a\in \partial M$, $O_{a, \epsilon}(1)$ stands for quantities uniformly bounded in $a$ and $\epsilon$. For \;$\epsilon$\; positive and small, $o_{\epsilon}(1)$ means quantities which tend to $0$ as $\epsilon$ tends to $0$. For $\l$ large and \;$a\in \partial M$, \;$O_{a, \l}(1)$ stands for quantities uniformly bounded in $a$ and $\l$. For \;$a\in \partial M$, \;$\epsilon$\; and \;$\delta$ positive and small, and $\l$ large, $O_{a, \epsilon, \delta}(1)$ and \;$O_{a, \l}(1)$\; stand respectively for quantities which are bounded uniformly in \;$a$, \;$\d$, and $\epsilon$, and in \;$a$\; and \;$\l$. For $a\in \partial M$, $\epsilon$ positive and small, and $\l$ large, $o_{a, \epsilon}(1)$ and $o_{a, \l}(1)$ stand respectively for quantities which tend to \;$0$\; uniformly in $a$ as $\epsilon$ tends to $0$, and as  $\l$ tends to $+\infty$. For \;$A\in (\partial M)^2$ and $\l$ large, $O_{A, \l}(1)$\; and $o_{A, \l}(1)$ stand respectively for quantities which are bounded uniformly in \;$A$\; and \;$\l$, and which tend to \;$0$ uniformly in \;$A$\; as \;$\l$\; tends to \;$+\infty$. For \;$p\in \N^*$, $A\in (\partial M)^p$, $\bar \alpha\in \D_{p-1}$, and \;$\l$ large, $O_{A, \bar \alpha, \l}(1)$ and $o_{A, \bar \alpha, \l}(1)$ stand respectively for quantities which are uniformly bounded in \;$p$, $A$, $\bar\alpha$, and $\l$ and for quantities which tend to \;$0$\; uniformly in \;$p$, $A$, and \;$\bar\alpha$ \;as \;$\l$\; tends to $+\infty$. For $x\in \R$, we will use the notation \;$O(x)$\; and \;$o(x)$\; to mean  respectively\;$|x|O(1)$\; and $|x|o(1)$ where \;$O(1)$\; and \;$o(1)$\; will be specified in all the contexts where they are used. Large positive constants are usually denoted by \;$C$\; and the value of \;$C$\; is allowed to vary from formula to formula and also within the same line. Similarly small positive constants are denoted by \;$c$\; and their values may vary from formula to formula and also within the same line.   The symbol \;$\sum_{i\neq j}$\;always means a double sum over the associated index set under the assumption \;$i\neq j$.                                                                                                                                                                                                                                                                                                                                                                                                                                                                                                                                                                                                                                                                                                                                                                                                                                                                                                                                                                                                                                                                                                    \vspace{6pt}

\noindent
For \;$X$ a topological space, \;$H_{*}(X)$\; will denote the singular homology of \;$X$ with \;$\Z_2$\; coefficients, and $H^*(X)$ for the cohomology. For \;$Y$ a subspace of \;$X$, $H_*(X, Y)$ will stand for the relative homology. The symbol \;$\frown$\; will denote the cap product between cohomology and homology. For a map $f:X\rightarrow Y$, with \;$X$\; and \;$Y$\; topological spaces, \;$f_*$\; stands for the induced map in homology, and \;$f^*$\; for the induced map in cohomology.
For $p\in \N$, we set
\begin{equation}\label{eq:defenergylevel}
W_p:=\{u\in W^{1, 2}_+(\ov M):\;\;\;\;\mathcal{E}_g(u)\leq (p+1)^{\frac{1}{n-1}}\mathcal{Q}(B^n)\}.
\end{equation}
where
\begin{equation}\label{eq:yamabehsphere}
\mathcal{Q}(B^n):=\mathcal{Q}(B^n, \partial B^n, g_{B^n})
\end{equation}
\vspace{6pt}

\noindent
For a Riemannian metric $\bar g$ defined on $\ov M$, we denote by \;$G_{\bar g}$\; the Green's function of $(L_{\bar g}, B_ {\bar g})$ satisfying the normalization
\begin{equation}\label{eq:greennorm}
\lim_{d_{\bar g}(a, x) \longrightarrow 0}(d_{\bar g}(a, x))^{n-2}G_{\bar g}(a, x)=1,
\end{equation} 
and set
\begin{equation}\label{eq:greenback}
G:=G_g.
\end{equation}
Using the existence of conformal normal Fermi coordinates (see \cite{marques}) and recalling that \;$\partial M$\;is totally geodesic in \;$(\ov M, g)$, we have that for every large positive integer \;$m$\; and for every \;$a \in \partial M$, there exists a positive function $u_a\in C^{\infty}(\ov M)$ such that the metric $g_a = u_a^{\frac{4}{n-2}}g$ verifies
\begin{equation}\label{eq:detga}
det g_a(x)=1 +O_{a, x}((d_{g_a}(x, a))^m)\;\;\text{for}\;\;\; x\in B^{g_a}_a( \varrho_a),
\end{equation}
with \;$O_{a, x}(1)$\; meaning bounded by a constant independent of \;$a$\; and \;$x$, \;$0<\varrho_a<\min\{\frac{inj_{g_a}(\ov M)}{10},\; \frac{inj_{\hat g_a}(\partial M)}{10}\}$. Moreover, we can take the family of functions $u_a$, $g_a$ and $\varrho_a$ such that
\begin{equation}\label{eq:varro0}
\text{the maps}\;\;\;a\longrightarrow u_a, \;g_a\;\;\text{are}\;\;C^0\;\;\;\text{and}\;\;\;\;\frac{1}{4}\geq \varrho_a\geq \varrho_0>0,
\end{equation}
for some small positive $\varrho_0$ satisfying $\varrho_0<\min\{\frac{inj_g(\ov M)}{10}, \;\frac{inj_{\hat g}(\partial M)}{10}\}$, and
\begin{equation}\label{eq:proua}
\begin{split}
&||u_a||_{C^2(\ov M)}=O_a(1),\;\;\frac{1}{\ov C^2} g\leq g_a\leq \ov C^2 g, \;\;\;a\in \ov M, \\\;&u_a(x)=1+ O_a(d^2_{g_a}(a, x))=1+O_a(d_{g}^2(a, x)) \;\;\text{for}\;\;x\in\;\;B_a^{g_a}(\varrho_0)\supset B_a(\frac{\varrho_0}{2\ov C}),\\&
u_a(a)=1,\;\text{and}\;\; H_{g_a}=0,
\end{split}
\end{equation}
for some large positive constant $\ov C$ independent of $a$, and for the meaning of $O_a(1)$ in \eqref{eq:proua}, see section \ref{eq:notpre}. For  \;$a\in \partial M$\; and \;$\epsilon$ positive, we recall that the standard bubbles of the geometric problem under study  are defined as follows
\begin{equation}\label{eq:standardbubbles}
\d_{a, \epsilon}(x):=\left(\frac{\epsilon}{(\epsilon+x_n)^2+|x^{'}|^2}\right)^{\frac{n-2}{2}}, \;\;x\in B_a^{g_a}(\varrho_0),
\end{equation}
where\; $(x^{'}, x_n)$\; is the Fermi normal coordinate of \;$x$\; with respect to \;$g_a$\; at $a$.
For $a\in \partial M$ and $0<r<\varrho_0$, we set also 
\begin{equation}\label{eq:greena}
G_a:=G_{g_a}, \;\;\;\;\exp_a^a=\exp_{a}^{g_a},\;\;\;\;\hat{\exp}_a^a=\exp_{a}^{\hat g_a},\;\;\;\;B_a^a(r):=B_a^{g_{a}}(r)\;\;\;\;\text{and}\;\;\;\hat B_a^a(r):=\hat B_a^{g_{a}}(r).             
\end{equation}
On the other hand, the conformal invariance properties of the couple conformal Laplacian and conformal Neumann operator imply
\begin{equation}\label{eq:invarpro}
\begin{split}
&\mathcal{E}_g(u)=\mathcal{E}_{g_{a}}(u_a^{-1}u), \;\;\;\;\oint_{\partial M}u^{\frac{2(n-1)}{n-2}}dS_g=\oint_{\partial M}(u_a^{-1}u)^{\frac{2(n-1)}{n-2}}dS_{g_a}\;\;\;\;\text{for}\;\;\;\;u\in W^{1, 2}_+(\ov M),\\ &G_g(x, y)=G_{a}(x, y)u_a(x)u_a(y),\;\; \;\;(x, y)\in \ov M^2\;\;\text{and}\;\;\;\;a\in \partial M.
\end{split}
\end{equation}
We also define the following quantities
\begin{equation}\label{eq:defc}
c_0:=(n-2), \;\;\;c_1:=\int_{\R^{n-1}}\left(\frac{1}{1+|x|^2}\right)^{n-1}dx,  \;\;\text{and}\;\;\;\;c_2:=4\frac{n-1}{n-2}\int_{\R^n_+}\left|\n\left[\left(\frac{1}{(1+x_n)^2+|x^{'}|^2}\right)^{\frac{n-2}{2}}\right]\right|^2dx.
\end{equation}
Furthermore, we set
\begin{equation}\label{defc3}
c_3:=\int_{\R^{n-1}}\left(\frac{1}{1+|x|^2}\right)^{\frac{n}{2}}dx,
\end{equation}
and define the following quantity which depends only on $(\ov M, g)$
\begin{equation}\label{eq:defcg}
c_g=\frac{c_3}{4c_1}\min_{((\partial M)^2)^*}G,
\end{equation}
and see above for the definition of \;$((\partial M)^2)^*$, \;$G$\; and \;$c_3$. We recall that the numbers $c_i$ ($i=0, 1, 2$) and \;$\mathcal{Q}(B^n)$\; verify the following relation 
\begin{equation}\label{eq:relationcy}
 c_2=c_0c_1\;\;\;\text{and}\;\;\;\mathcal{Q}(B^n)=\frac{c_2}{c_1^{\frac{n-2}{n-1}}}.
\end{equation}
Moreover, for \;$p\in \N^*$, $A:=(a_1, \cdots, a_p)\in (\partial M)^p$, \;$\bar \l:=(\l_1, \cdots, \l_p)\in (\R_+)^p$, we associate the following quantities (which appear in the analysis of diverging PS sequences of the Euler-Lagrange functional \;$\mathcal{E}_g$)
 \begin{equation}\label{eq:varepsilonij}
     \varepsilon_{i, j}:=\varepsilon_{i, j}(A, \bar \l):=\frac{c_{3}}{(\frac{\l_i}{\l_j}+\frac{\l_j}{\l_i}+\lambda_i\l_jG^{\frac{2}{2-n}}(a_{i},a_{i}))^{\frac{n-2}{2}}}, \;\;i, j=1, \cdots, p, \;\;i\neq j.                                                                                                                                                                                                                                                                                                                                                                                                                                                                                                                                                                                                                                                                                                                                                                                                                                                                                                                                                                                                                                                                                                                                                                                           \end{equation}
                                                                  
\vspace{6pt}

\noindent
Now, we are going to present some topological properties of the space of formal barycenter of $\partial M $ that we will need for our algebraic topological argument for existence. To do that,  for \;$p\in \N^*$, we recall that the set of formal barycenters of \;$\partial M$\; of order $p$ is defined as follows
 \begin{equation}\label{eq:barytop}
B_{p}(\partial M):=\{\sum_{i=1}^{p}\alpha_i\d_{a_i}:\:\;\;\;a_i\in \partial M, \;\alpha_i\geq 0,\;\; i=1,\cdots, p,\;\,\sum_{i=1}^{p}\alpha_i=1\}, 
\end{equation}
and set
\begin{equation}
B_0(\partial M)=\emptyset. 
\end{equation}
Furthermore, we have the existence of  \;$\Z_2$\; orientation classes \;$w_p\in H_{np-1}(B_{p}(\partial M), B_{p-1}(\partial M))$\; and that the cap product acts as follows 
\begin{equation}\label{eq:actioncap}
\begin{CD}
 H^l(B_p(\partial M)\setminus B_{p-1}(\partial M))\times H_k(B_{p}(\partial M), B_{p-1}(\partial M))@>\frown>> H_{k-l}(B_{p}(\partial M), B_{p-1}(\partial M)).
 \end{CD}
\end{equation}
Moreover, there holds
\begin{equation}\label{eq:purem}
B_p(\partial M)\setminus B_{p-1}(\partial M)\simeq ((\partial M)^p)^*\times_{\sigma_p}\D_{p-1}^*,
\end{equation}
where
$$
((\partial M)^p)^*:=(\partial M)^p)\setminus F((\partial M)^p),
$$
and
$$
\D_{p-1}^*:=\{(\alpha_1, \cdots, \alpha_p)\in \D_{p-1}:\;\; \;\alpha_i>0, \;\;\forall\;i=1, \cdots, p\}.
$$
On the other hand, since $\partial M$ is a closed $(n-1)$-dimensional manifold, then we have 
\begin{equation}\label{eq:defom}
\text{an orientation class}\,\;0\neq O^{*}_{\partial M}\in H^{n-1}(\partial M).
\end{equation}
Furthermore, there is a natural way to inject $\partial M$ into $((\partial M)^p)^*\times_{\sigma_p}\D_{p-1}^*$, namely an injection
\begin{equation}\label{eq:inj}
i: \partial M\longrightarrow ((\partial M)^p)^*\times_{\sigma_p}\D_{p-1}^*,
\end{equation}
such that
\begin{equation}\label{eq:defom1}
i^*(O^*_p)=O^*_{\partial M}\;\; \text{with}\;\;0\neq O^{*}_p\in H^{n-1}(((\partial M)^p)^*\times_{\sigma_p}\D_{p-1}^*).
\end{equation}
Identifying $O^*_{\partial M}$ and $O^*_p$ via  \eqref{eq:defom1}, and using \eqref{eq:actioncap} and \eqref{eq:purem}, we have the following well-know formula, see \cite{kk}.
\begin{lem}\label{eq:transfert}
There holds 
$$
\begin{CD}
 H^{n-1}(((\partial M)^p)^*\times_{\sigma_p}\D_{p-1}^*)\times H_{np-1}(B_{p}(\partial M), B_{p-1}(\partial M))@>\frown>> H_{np-n}(B_{p}(\partial M), B_{p-1}(\partial M))\\@>\partial>>H_{np-n-1}(B_{p-1}(\partial M), B_{p-2}(\partial M)),
 \end{CD}
$$
and 
\begin{equation}\label{eq:classt}
\omega_{p-1}=\partial(O^*_{\partial M}\frown w_p).
\end{equation}
\end{lem}
\vspace{8pt}

\noindent
Next, we are going to discuss some important properties of the Chen\cite{chen}'s bubbles. Using the techniques of Brendle\cite{bre2}, Chen\cite{chen} has introduced a family of bubbles which verify the same properties as the Brendle\cite{bre2}'s bubbles and the Brendle-Chen\cite{bs}'s bubbles. Indeed, for $\delta$ small, she defines a family of bubbles $v_{a, \epsilon, \delta}$ (see page 16 in \cite{chen}), $a\in \partial M$\; and $\epsilon$\; positive and small such that they can replace the standard bubbles in the analysis of diverging PS sequences of $\mathcal{E}_g$ and more importantly verify a sharp energy estimate. Precisely, $v_{a, \epsilon, \delta}$\; is defined as a suitable perturbation of the standard bubbles glued with an appropriate scale of the Green's function \;$G_a$\; centered at $a$ as follows
\begin{equation}\label{eq:chenbubbles}
\begin{split}
v_{a, \epsilon, \delta}(\cdot)
=
\chi_{\delta}(\cdot)(\d_{a, \epsilon}(\cdot)+w_{a, \epsilon}(\cdot))
+
(1-\chi_{\delta}(\cdot))\epsilon^{\frac{n-2}{2}}G_{a}(a, \cdot),
\end{split}\end{equation}
where 
\begin{equation}\label{eq:cutoff}
\chi_{\delta}(x):=\chi\left(\frac{d_{g_a}(a, x)}{\delta}\right),
\end{equation}
and\;$\chi$ is a cut-off function defined on $\bar \R_+$ satisfying $\chi$ is non-negative, $\chi(t)=1$ if $t\leq 1$ and $\chi(t)=0$ if $t\geq 2$, $\delta_{a, \epsilon}$ is defined as in \eqref{eq:standardbubbles}, $G_{a}(a, \cdot)$ is defined as in \eqref{eq:greena}, and in Fermi normal coordinates around $a$ with respect to $g_a$, we have that $w_{a, \epsilon}$ satisfies the following pointwise estimate
\begin{equation}\label{eq:wdecay}
\vert \partial^{\beta}w_{a, \epsilon}(x)\vert
\leq C_n(|\beta|)
\frac
{\epsilon^{\frac{n-2}{2}}}
{(\epsilon^{2}+r^{2})^{\frac{n-4+\beta}{2}}}\;\;\text{with}\;\;\;r=d_{g_a}(a, x)\; \;\text{and}\; \;\;x\in B_a^a(\varrho_0),
\end{equation}
where $\varrho_0$ is as in \eqref{eq:proua} and $C_n(|\beta|)$ is a large positive constant which depends only on $n$ and $|\beta|$. Furthermore, $v_{a, \epsilon, \delta}$\;verifies the following energy estimate which is a weak form of Proposition 9 in  \cite{chen}, but sufficient for the purpose of this paper.
\begin{lem}\label{eq:brenchenenergy}
There exists \;$0<\d_0\leq\varrho_0$ small such that for every \;$0<2\epsilon\leq\delta\leq \d_0$ and for every $a\in \partial M$, there holds
\begin{equation}\label{eq:brenchenenergyest}
\begin{split}
\int_M\left(4\frac{n-1}{n-2}|\n_{g_a} v_{a, \epsilon, \delta_0}|^2+R_{g_a}v_{a, \epsilon, \delta}^2dV_{g_a}\right)\leq \mathcal{Q}(B^n)\left(\oint_{\partial M}v_{a, \epsilon, \delta}^{\frac{2(n-1)}{n-2}}dS_{g_a}\right)^{\frac{n-2}{n-1}}-\epsilon^{n-2}\mathcal{I}(a, \delta)\\+O_{a, \epsilon, \delta}(\d^2\epsilon^{n-2}+\epsilon^{n-1}\delta^{-n+1}),
\end{split}
\end{equation}
where \;$\mathcal{Q}(B^n)$\; is defined by \eqref{eq:yamabehsphere}, \;$\mathcal{I}(a, \delta)$\; is a flux integral verifying \;$\mathcal{I}(a, \delta)=O_{a, \delta}\left(1\right)$, and for the meaning of $O_{a, \delta}\left(1\right)$ and $O_{a, \epsilon, \delta}\left(1\right)$, see Section \ref{eq:notpre}.
\end{lem}
\vspace{6pt}

\noindent
On the other hand, using the work of Almaraz\cite{almaraz2} (Lemma 3.15 in \cite{almaraz2}), we have that the\;$v_{a, \epsilon, \delta}$'s\; verify the following interaction estimates.
\begin{lem}\label{eq:brencheninteract}
There exists a large constant $C_1>0$ such that for every \;$2\epsilon_1\leq 2\epsilon_2\leq\delta^2\leq \delta_0$\; and every $a_1, a_2\in \partial M$, there holds
 \begin{equation}\label{eq:brencheninteractest}
 \begin{split}
 &\int_M v_{a_1,  \epsilon_1, \delta}\left|-4\frac{n-1}{n-2}\D_{g_{a_2}} v_{a_2, \epsilon_2, \delta}+R_{g_{a_2}}v^2_{a_2, \epsilon_2, \delta}\right|dV_{g_{a_2}}+\oint_{\partial M}v_{a_1, \epsilon_1, \delta}\left|4\frac{n-1}{n-2}\frac{\partial v_{a_2, \epsilon_2, \delta}}{\partial n_{g_{a_2}}}-c_0v_{a_2, \epsilon_2,\delta}^{\frac{n}{n-2}}\right|dS_{g_{a_2}}\\&\leq C_1\left(\delta+\frac{\epsilon_2}{\delta}\right)\left(\frac{\epsilon_2^2+d_{g_{a_2}}^2(a_1, a_2)}{\epsilon_1\epsilon_2}\right)^{\frac{2-n}{2}},
 \end{split}
\end{equation}
where \;$c_0$\; is defined by \eqref{eq:defc}.
\end{lem}
\vspace{6pt}

\noindent
Furthermore, using \eqref{eq:chenbubbles}, it is easy to see that the following estimate holds.
\begin{lem}\label{eq:bsvolest}
Assuming that \;$0<\epsilon\leq\delta_0^{\frac{n-1}{2}}$, and \;$a\in \partial M$, then we have
\begin{equation}\label{eq:estconformalvol}
\oint_{\partial M}v_{a, \epsilon, \epsilon^{\frac{2}{n-1}}}^{\frac{2(n-1)}{n-2}}dS_{g_a}=c_1+o_{a, \epsilon}(1),
\end{equation}
where \;$c_1$\; is as in \eqref{eq:defc}, and for the meaning of \;$o_{a, \epsilon}(1)$, see Section \ref{eq:notpre}.
\end{lem}
\vspace{6pt}

\noindent
Thus, setting 
\begin{equation}\label{eq:val}
 v_{a}^{ \l}:=v_{a, \frac{1}{\l}, (\frac{1}{\l})^{\frac{2}{n-1}}}, \;\;a\in \partial M,\;\;\;\;\l\geq \frac{2}{\delta_0^{\frac{n-1}{2}}}
\end{equation}
and
\begin{equation}\label{eq:varphialb}
\varphi_{a, \l}:=u_av_{a}^{\l}, \;\;a\in \partial M, \;\;\;\;\l\geq \frac{2}{\delta_0^{\frac{n-1}{2}}},
\end{equation}
where \;$u_a$\; is as in \eqref{eq:detga} and \;$\delta_0$ is still given by Lemma \ref{eq:brenchenenergy}, we have clearly that Lemma \ref{eq:brenchenenergy}, Lemma \ref{eq:brencheninteract}, and Lemma \ref{eq:bsvolest} combined with \eqref{eq:invarpro} imply the following Lemmata which will play an important role in our application of the barycenter technique of Bahri-Coron\cite{bc}.
\begin{lem}\label{eq:bubbleestl}
Assuming that $a\in \partial M$ and \;$\l\geq \frac{2}{\delta_0^{\frac{n-1}{2}}}$, then the following estimate holds
\begin{equation}\label{eq:bubblesestenergy}
\mathcal{E}_g(\varphi_{a, \l})\leq \mathcal{Q}(B^n)\left(1+O_{a, \l}\left(\frac{1}{\l^{n-2}}\right)\right),
\end{equation}
where $\mathcal{Q}(B^n)$ is as in \eqref{eq:yamabehsphere} and for the meaning of \;$O_{a, \l}\left(1\right)$, see Section \ref{eq:notpre}.
\end{lem}
\begin{lem}\label{eq:bubbleinteractl}
 There exits a large constant $C_2>0$ such that for every $a_1, a_2\in \partial M$, and for every $\l\geq \frac{2}{\delta_0^{\frac{n-1}{2}}}$, we have
 \begin{equation}\label{eq:bubbleinteractest}
 \begin{split}
 &\int_M \varphi_{a_1,  \l}\left|-4\frac{n-1}{n-2}\D_{g_{a_2}} \varphi_{a_2, \l}+R_{g_{a_2}}\varphi^2_{a_2, \l}\right|dV_{g_{a_2}}+\oint_{\partial M}\varphi_{a_1, \l}\left|4\frac{n-1}{n-2}\frac{\partial \varphi_{a_2, \l}}{\partial n_{g_{a_2}}}-c_0\varphi_{a_2, \l}^{\frac{n}{n-2}}\right|dS_{g_{a_2}}\\&\leq C_1\left[\left(\frac{1}{\l}\right)^{\frac{2}{n-1}}+\left(\frac{1}{\l}\right)^{\frac{n-3}{n-1}}\right]\left(1+\l^2d_{g_{a_2}}^2(a_1, a_2)\right)^{\frac{2-n}{2}},
 \end{split}
 \end{equation}
 where \;$c_0$\; is as in \eqref{eq:defc}.
\end{lem}\begin{lem}\label{eq:bsvolestl}
Assuming that  $a\in \partial M$ and $\l\geq \frac{1}{\delta_0^{\frac{n-1}{2}}}$, then there holds
\begin{equation}\label{eq:estconformalvol}
\oint_{\partial M} \varphi_{a, \l}^{\frac{2(n-1)}{n-2}} dS_{g_a}=c_1+o_{a, \l}(1),
\end{equation}
where \;$c_1$\; is as in \eqref{eq:defc} and for the meaning of \;$o_{a, \l}(1)$, see Section \ref{eq:notpre}.
\end{lem}
\vspace{6pt}

\noindent
On the other hand, using \eqref{eq:chenbubbles}-\eqref{eq:wdecay}, and \eqref{eq:val}, we have that \;$v_{a}^{\l}$\; decomposes as follows 
\begin{equation}\label{eq:decompval}
v_{a}^{ \l}(\cdot)=\chi^{\l}(\cdot)\left(\delta_{a}^{ \l}(\cdot)+w_{a}^{\l}(\cdot)\right)\left(1-\chi_{\l}(\cdot)\right)\frac{G_{a}(a, \cdot)}{\l^\frac{{n-2}}{2}}
\end{equation}
where
\begin{equation}\label{eq:wal}
w_{a}^{ \l}:=w_{a, \frac{1}{\l}},\;\;\;\delta_{a}^{\l}:=\delta_{a, \frac{1}{\l}}, \;\;\text{and}\;\;\;\chi^{\l}=\chi_{(\frac {1}{\l})^{\frac{2}{n-1}}},
\end{equation}
and \;$w_{a}^{\l}$\; satisfies the following pointwise estimate
\begin{equation}\label{eq:walest}
\vert \partial^{\beta}w_{a}^{\l}(x)\vert
\leq C_n(|\beta|)
\frac
{\l^{\frac{n-6}{2}+\beta}}
{(1+\l^2r^{2})^{\frac{n-4+\beta}{2}}}\;\;\text{with}\;\;\;r=d_{g_a}(a, x)\; \;\text{and}\; \;\;x\in B_a^a(\varrho_0).
\end{equation}
Now, for $p\in \N^*$,  and $A:=(a_1, \cdots, a_p)\in (\partial M)^p$ and $\l\geq \frac{2}{\delta_0^{\frac{n-1}{2}}}$, we associate the following quantities
\begin{equation}\label{eq:intweight}
\epsilon_{i, j}:=\epsilon_{i, j}(A, \l):= \oint_{\partial M} \varphi_{a_i, \l}^{\frac{n}{n-2}}\varphi_{a_j, \l}dS_g, \;\;\;i, j=1, \cdots, p, \;\;i\neq j.
\end{equation}
Using \eqref{eq:varphialb}, \eqref{eq:decompval}-\eqref{eq:intweight}, we have the following Lemma which provides self and interaction estimates, and a relation between \;$\epsilon_{i, j}(A, \l)$\; and \;$\varepsilon_{i, j}(A, \bar \l)$ with $\bar \l:=(\l, \cdots, \l)$, and for the meaning of \;$\varepsilon_{i, j}(A, \bar \l)$ see \eqref{eq:varepsilonij}.
\begin{lem}\label{eq:interactestl}
Assuming that $p \in \N^*$, $A:=(a_1,\cdots, a_p)\in (\partial M)^p$ and $\l\geq \frac{2}{\d_0^{\frac{n-1}{2}}}$, then \\
1) For every $i, j=1, \cdots, p$ with $i\neq j$, we have \\
i) 
\begin{equation*}
\begin{split}
\epsilon_{i,j}\longrightarrow 0
\Longleftrightarrow
\varepsilon_{i,j}\longrightarrow 0,
\end{split}
\end{equation*}
where $\varepsilon_{i, j}:=\varepsilon_{i, j}(A, \bar  \l)$ with $\bar \l:=(\l, \cdots, \l)$ and $\epsilon_{i, j}:=\epsilon_{i, j}(A, \l)$, and for their definitions see respectively \eqref{eq:varepsilonij} and \eqref{eq:intweight}.\\
ii)\\
There exists $\;0<C_3<\infty$\; independent of \;$p$, \;$A$ and \;$\l$ such that the following estimate holds
\begin{equation}
\begin{split}
C_3^{-1}<\frac{\epsilon_{i,j}}{\varepsilon_{i,j}}<C_3.
\end{split}
\end{equation} 
iii)\\
If \;$\varepsilon_{i,j} \longrightarrow 0$, then 
\begin{equation*}
\begin{split}
\epsilon_{i,j}=(1+o_{ \varepsilon_{i, j}}(1))\varepsilon_{i,j},
\end{split}
\end{equation*}
and for the meaning of $o_{\varepsilon_{i, j}}(1)$, see Section \ref{eq:notpre}.\\
2) For every \;$i=1, \cdots, p$, there holds
\begin{equation*}
\begin{split}
\langle L_{g}\varphi_{a_i, \l}, \varphi_{a_i, \l}\rangle+\langle B_g \varphi_{a_i, \l}, \varphi_{a_i, \l}\rangle=c_0(1+o_{a
_i, \l}(1))\oint_{\partial M} \varphi^{\frac{2(n-1)}{n-2}}_{a_i, \l}dS_g,
\end{split}
\end{equation*}
where $c_0$ is given by \eqref{eq:defc} and for the meaning of $o_{a_i, \l}(1)$, see Section \ref{eq:notpre}.\\
3) For every \;$i, j=1,\cdots, p$ with $i\neq j$, there holds
\begin{equation*}
\begin{split}
\langle L_{g}\varphi_{a_i}, \varphi_{a_j}\rangle+\langle B_g\varphi_{a_i, \l}, \varphi_{a_j, \l}\rangle=(1+o_{A_{i, j}, \l}(1))c_{0}\epsilon_{i,j},\;\;\;\;\text{and}\;\;\;\;\;
\epsilon_{j,i}=(1+o_{A_{i, j}, \l}(1))\epsilon_{i,j},
\end{split}
\end{equation*}
where $A_{i, j}:=(a_i, a_j)$ and for the meaning of $o_{A_{i, j}, \l}(1)$, see Section \ref{eq:notpre}.
\end{lem}
\begin{pf}
To prove Lemma \ref{eq:interactestl}, we use the same strategy as in \cite{martndia1}. First of all, to simplify notation, for every  $i=1, \cdots, p$, we set
 \begin{equation}\label{eq:auxilary}
 \begin{split}
&\varphi_{i}:=\varphi_{a_i, \l}, \;\; v_i:=v_{a_i}^{ \l}\;\;\d_i:=\delta_{a_i}^{\l}, \;\;w_{i}:=w_{a_i}^{ \l}, \;\;G_i(\cdot):=G_{a_i}(a_i, \cdot),\\& B_i:=\hat B_{a_i}^{a_i}(\d_0), \;\;\;\exp_i:=\hat{\exp}_{a_i}^{a_i}\; \;\text{and}\;\;u_i:=u_{a_i},
\end{split}
\end{equation}
and for the meaning of $\hat B_{a_i}^{a_i}(\d_0)$ and $\hat{\exp}_{a_i}^{a_i}$, see Section \ref{eq:notpre}.
Next, using  \eqref{eq:walest} and \eqref{eq:auxilary}, we have that $v_i$ verifies the following pointwise estimate
\begin{equation}\begin{split}
v_{i}
= &
\left(1+o_{a_i, \l}(1)\right)\left(\chi^{\l}\d_i+(1-\chi^{\l})\frac{G_i}{\l^{\frac{n-2}{2}}}\right).
\end{split}\end{equation}
Moreover, using \eqref{eq:greennorm}, \eqref{eq:standardbubbles}, \eqref{eq:wal}, and \eqref{eq:auxilary}, we obtain (on \;$ \partial M$)
\begin{equation}\label{eq:estp1}
\begin{split}
\chi^{\l}\d_{i}
= &
\chi^{\l}\left(\frac{\l}{1+\l^2r^{2}}\right)^{\frac{n-2}{2}}
=
\chi^{\l}\left(\frac{\l}{1+\l^2G_{i}^{\frac{2}{2-n}}\frac{r^{2}}{G_{i}^{\frac{2}{2-n}}}}\right)^{\frac{n-2}{2}} \\
= &
\left(1+o_{a_i, \l}(1)\right)\chi^{\l}\left(\frac{\l}{1+\l^2G_{i}^{\frac{2}{2-n}}}\right)^{\frac{n-2}{2}}
\end{split}\end{equation}
where $r$ is as in \eqref{eq:wdecay} with $a$ replaced by $a_i$, and
\begin{equation}\label{eq:estp2}
\begin{split}
\left(1-\chi^{\l}\right)\left(\frac{\l}{1+\l^2G_{i}^{\frac{2}{2-n}}}\right)^{\frac{n-2}{2}}
= &
\left(1-\chi^{\l}\right)\frac{G_{i}}{\l^{\frac{n-2}{2}}}\left(\frac{1}{1+O_{a_i, \l}(\vert \frac{1}{\l}\vert^{\frac{2(n-3)}{n-1}})}\right) \\
= &
\left(1+o_{a_i, \l}(1)\right)\left(1-\chi^{\l}\right)\frac{G_i}{\l^{\frac{n-2}{2}}}.
\end{split}\end{equation}
Hence, combining \eqref{eq:estp1} and \eqref{eq:estp2}, we obtain
\begin{equation}\label{eq:viest}
\begin{split}
v_{i}=\left(1+o_{a_i, \l}(1)\right)\left(\frac{\l}{1+\l^2G_{i}^{\frac{2}{2-n}}}\right)^{\frac{n-2}{2}}.
\end{split}\end{equation}
Now, using  \eqref{eq:proua},  \eqref{eq:varphialb}, \eqref{eq:intweight}, \eqref{eq:auxilary}, and \eqref{eq:viest}, we derive the following estimate for $\epsilon_{i, j}$ ($i, j=1, \cdots, p$ and $i\neq j$)
\begin{equation}\label{eq:41}
\begin{split}
\epsilon_{i,j}
= &
\oint_{\partial M} v_{i}^{\frac{n}{n-2}}v_{j}\frac{u_{j}}{u_{i}}dS_{g_{a_{i}}} \\
= &
\int_{B_i}v_{i}^{\frac{n}{n-2}}v_{j}\frac{u_{j}}{u_{i}}dS_{g_{a_{i}}}
+
O_{A_{i, j}, \l}\left(\left(\frac{1}{\l}\right)^{n-1}\right)\\=&\left(1+o_{A_{i, j}, \l}(1)\right)
u_{j}(a_{i})\int_{B^{\l}(0)}\left(\frac{1}{1+|x|^{2}}\right)^{\frac{n}{2}}
\left[\frac{1}{1+\l^{2}G_{j}^{\frac{2}{2-n}}(\exp_{i}(\frac{x}{\l}))}\right]^{\frac{n-2}{2}} dx\\&+
O_{A_{i, j}, \l}\left(\left(\frac{1}{\l}\right)^{n-1}\right),
\end{split}\end{equation}
where 
\begin{equation}
B^{\l}(0):=B^{\R^{n-1}}_0(\l\d_0),      
\end{equation}
and for the meaning of $B^{\R^{n-1}}_0(\l\d_0)$, see Section \ref{eq:notpre}. From \eqref{eq:41} it follows that 
\begin{equation}\label{eq:equivalence0}
\begin{split}
\epsilon_{i,j}\longrightarrow 0
\Longleftrightarrow
\l^2 G_{j}^{\frac{2}{2-n}}(a_{i})\longrightarrow+ \infty
\Longleftrightarrow
\varepsilon_{i,j}\longrightarrow 0.
\end{split}\end{equation}
Thus, we have that the proof of i) of point 1) is complete. Now, since $\epsilon_{i,j}$ and $\varepsilon_{i,j}$ are bounded by definition, then thanks to \eqref{eq:equivalence0}, to prove ii) of point 1), we can assume without loss of generality that 
\begin{equation}\label{eq:large1}
\begin{split}
\l^{2}G_{j}^{\frac{2}{2-n}}(a_{j})\gg 1. 
\end{split}
\end{equation} 
Thus under the latter assumption, setting
\begin{equation}
\mathcal{A}=B^{\R^{n-1}}_{0}(\gamma\l\sqrt{G_{j}^{\frac{2}{2-n}}(a_{i})})
\end{equation} 
for $\gamma>0$ small and using Taylor expansion, we obtain that the following estimate holds on\; $\mathcal{A}$
\begin{equation}\label{eq:45}
\begin{split}
&\left(\frac{1}{1+\l^{2}G_{j}^{\frac{2}{2-n}}(\exp_{i}(\frac{x}{\l}))}\right)^{\frac{n-2}{2}}=  
\left(
\frac
{1}
{1+\l^{2}G_{j}^{\frac{2}{2-n}}(a_{i})}
\right)^{\frac{n-2}{2}}
\left(
\frac
{1}
{1
+
\frac
{\l^{2}\left[G_{j}^{\frac{2}{2-n}}(\exp_{i}(\frac{x}{\l}))-G_{j}^{\frac{2}{2-n}}(a_{i})\right]}{1+\l^{2}G_{j}^{\frac{2}{2-n}}(a_{i})}}
\right)^{\frac{n-2}{2}} \\
& = 
\left(
\frac
{1}
{1+\l^{2}G_{j}^{\frac{2}{2-n}}(a_{i})}
\right)^{\frac{n-2}{2}}
\left(
\frac
{1}
{1
+
\frac
{\l\nabla G_{j}^{\frac{2}{2-n}}(a_{i})x+O(\vert x \vert^{2})}{1+\l^{2}G_{j}^{\frac{2}{2-n}}(a_{i})}}
\right)^{\frac{n-2}{2}} \\
& = 
\left(
\frac
{1}
{1+\l^{2}G_{j}^{\frac{2}{2-n}}(a_{i})}
\right)^{\frac{n-2}{2}}-\frac{(n-2)}{2}\l
\left(
\frac
{1}
{1+\l^{2}G_{j}^{\frac{2}{2-n}}(a_{i})}
\right)^{\frac{n}{2}}
\nabla G_{j}^{\frac{2}{2-n}}(a_{i})x
+
O\left(\frac{\vert x \vert^{2}}{(1+\l^{2}G_{j}^{\frac{2}{2-n}}(a_{i}))^{\frac{n}{2}}}\right).
\end{split}\end{equation}
Now, combining  \eqref{eq:41} and \eqref{eq:45}, we obtain
\begin{equation}\label{eq:47}
\begin{split}
\epsilon_{i,j}
= &
(1+o_{A_{i, j}, \l}(1))\frac{u_{j}(a_{i})c_3}{\left(1+\l^{2}G_{j}^{\frac{2}{2-n}}(a_{i})\right)^{\frac{n-2}{2}}}
+
o_{\varepsilon_{i, j}}(\varepsilon_{i,j}) \\
& +
\int_{\mathcal{A}^{c}\cap B^{\l}}
\left(\frac{1}{1+|x|^{2}}\right)^{\frac{n}{2}}
\left[\left(\frac{1}{1+\l^{2}G_{j}^{\frac{2}{2-n}}}\right)^{\frac{n-2}{2}}\right]\circ \exp_{i}(\frac{x}{\l})\;dx.
\end{split}\end{equation}
Next, using \eqref{eq:invarpro}, \eqref{eq:large1}, and Taylor expansion, we derive that
\begin{equation}\label{eq:48}
\begin{split}
\frac{u_{j}(a_{i})c_3}{\left(1+\l^{2}G_{j}^{\frac{2}{2-n}}(a_{i})\right)^{\frac{n-2}{2}}}
= &
c_3\left(1+o_{\varepsilon_{i, j}}(1)\right)u_{j}(a_{i})\frac{G_{j}(a_{i})}{\l ^{n-2}}\\
= &
c_3\left(1+o_{\varepsilon_{i, j}}(1)\right)\frac{G(a_{i},a_{j})}{\l^{n-2}}
=
\left(1+o_{\varepsilon_{i, j}}(1)\right)\varepsilon_{i,j}.
\end{split}
\end{equation}
Thus, combining \eqref{eq:equivalence0}, \eqref{eq:large1},  \eqref{eq:47}, and \eqref{eq:48}, we obtain
\begin{equation}\begin{split}
\epsilon_{i,j}
= &
\left(1+o_{\varepsilon_{i, j}}(1)\right)\varepsilon_{i,j} 
+
I_{\mathcal {A}^{c}},
\end{split}\end{equation}
where
\begin{equation}\begin{split}
I_{\mathcal{A}^c}=&
\int_{\mathcal{A}^{c}\cap B^{\l}(0)}
\left(\frac{1}{1+|x|^{2}}\right)^{\frac{n}{2}}
\left[\left(\frac{1}{1+\l^{2}G_{j}^{\frac{2}{2-n}}}\right)^{\frac{n-2}{2}}\right]\circ \exp_{i}(\frac{x}{\l})\;dx,
\end{split}\end{equation}
and \;$\mathcal{A}^c=\R^{n-1}\setminus \mathcal{A}$. Hence, to end the proof of ii) of point 1) and to prove iii) of point 1), we are going to show that $
I_{\mathcal{A}^c}$ satisfies
\begin{equation}\label{eq:iaest}
I_{\mathcal{A}^c}=o_{\varepsilon_{i, j}}(\varepsilon_{i, j}).
\end{equation}
In order to do that, we first decompose \;$\mathcal{A}^{c}$\; into
\begin{equation}\label{eq:defb}
\begin{split}
\mathcal{B}
=
\{x\in \R^{n-1}:\;\;
\gamma\l\sqrt{G_{j}^{\frac{2}{2-n}}(a_{i})}\leq \vert x \vert 
\leq 
\gamma^{-1}\l\sqrt{G_{j}^{\frac{2}{2-n}}(a_{i})}
\}
\end{split}\end{equation}
and
\begin{equation}\label{eq:deftildec}
\begin{split}
\mathcal{C}
=
\{x\in \R^{n-1}:\;\;
\vert x \vert > \gamma^{-1}\l\sqrt{G_{j}^{\frac{2}{2-n}}(a_{i})}
\},
\end{split}\end{equation}
and have
\begin{equation}\begin{split}
I_{\mathcal {A}^c}
=
I_{\mathcal{B}}
+
I_{\mathcal{C}},
\end{split}\end{equation}
where
\begin{equation}\label{eq:ib}
I_{\mathcal{B}}:=\int_{\mathcal{B}\cap B^{\l}(0)}
\left(\frac{1}{1+|x|^{2}}\right)^{\frac{n}{2}}
\left[\left(\frac{1}{1+\l^{2}G_{j}^{\frac{2}{2-n}}}\right)^{\frac{n-2}{2}}\right]\circ \exp_{i}(\frac{x}{\l})\;dx
\end{equation}
and
\begin{equation}\label{eq:ic}
I_{\mathcal{C}}:=\int_{\mathcal{C}\cap B^{\l}(0)}
\left(\frac{1}{1+|x|^{2}}\right)^{\frac{n}{2}}
\left[\left(\frac{1}{1+\l^{2}G_{j}^{\frac{2}{2-n}}}\right)^{\frac{n-2}{2}}\right]\circ \exp_{i}(\frac{x}{\l})\;dx.
\end{equation}
To prove \eqref{eq:iaest}, we are going to estimate separately $I_{\mathcal{B}}$ and $I_{\mathcal{C}}$. We start with $I_{\mathcal{B}}$. Using \eqref{eq:defb} and \eqref{eq:ib}, we have clearly that \;$I_{\mathcal{B}}$\; verifies the following estimate
\begin{equation}\label{eq:ibest1}
\begin{split}
I_{\mathcal{B}}
\leq 
\frac{C_{\gamma}}{\left(1+\l^{2}G_{j}^{\frac{2}{2-n}}(a_{i})\right)^{\frac{n}{2}}} \times
\int_{\mathcal{B}\cap B^{\l}}
\left[\left(\frac{1}{1+\l^{2}G_{j}^{\frac{2}{2-n}}}\right)^{\frac{n-2}{2}}\right]\circ \exp_{i}(\frac{x}{\l})\:dx,
\end{split}
\end{equation}
for some large positive constant \;$C_{\gamma}$ depending only on \;$\gamma$. Thus, rescaling and changing coordinates via $\exp_{j}\circ \exp_{i}^{-1}$ (if necessary), we have that \eqref{eq:ibest1} implies
\begin{equation}\label{eq:ibest2}
\begin{split}
I_{\mathcal{B}}
\leq & 
\hat C_{\gamma}\varepsilon_{i,j}^{\frac{n}{n-2}}
\int_{[\vert x \vert \leq \l\tilde C_{\gamma}d_{g}(a_{i},a_{j})]}
\left(\frac{1}{1+|x|^{2}}\right)^{\frac{n-2}{2}}dx
\leq
\bar C_{\gamma}\varepsilon_{i,j}^{\frac{n-1}{n-2}},
\end{split}\end{equation}
for some large positive constants \;$\hat C_{\gamma}$, $\tilde C_{\gamma}$, and \;$\bar C_{\gamma}$\; which are depending only on \;$\gamma$.
Finally, we estimate \;$I_{\mathcal{C}}$. To do that, we fix \;$\gamma>0$ \;sufficiently small and use \eqref{eq:equivalence0}, \eqref{eq:large1}, and \eqref{eq:ic} to obtain
\begin{equation}\label{eq:icest1}
\begin{split}
I_{\mathcal{C}}
\leq 
\frac{C_{\gamma}}{\left(1+\l^{2}G_{j}^{\frac{2}{2-n}}(a_{i})\right)^{\frac{n-2}{2}}}
\int_{\mathcal{C}}
\left(\frac{1}{1+|x|^{2}}\right)^{\frac{n}{2}}dx
=
o_{\varepsilon_{i, j}}(\varepsilon_{i,j}),
\end{split}\end{equation}
for some large constant \;$C_{\gamma}$\; depending only on \;$\gamma$. Hence \eqref{eq:ibest2} and \eqref{eq:icest1} imply \eqref{eq:iaest}, thereby ending the proof of point 1). On the other hand, we have clearly that point 2) follows from Lemma \ref{eq:bubbleinteractl} and Lemma \ref{eq:bsvolestl}. Furthermore, the first equation of point 3) follows from Lemma \ref{eq:bubbleinteractl}, and ii) of point 1), while the second equation follows from the first equation and from the self-adjointness of \;$(L_g, B_g)$.
\end{pf}
\vspace{6pt}

\noindent
Now, using \eqref{eq:standardbubbles}, \eqref{eq:varphialb}, \eqref{eq:decompval}-\eqref{eq:walest}, we have the following interaction type estimate.
\begin{lem}\label{eq:alphabetainteract}
Assuming that $p\in \N^*$, $A:=(a_1, \cdots, a_p)\in (\partial M)^p$ and $\l\geq \frac{2}{\delta_0^{\frac{n-1}{2}}}$, then for every $i, j=1, \cdots, p$ with $i\neq j$, there holds
\begin{equation}\label{eq:betaint}
\begin{split} 
\oint_{\partial M} \varphi_{a_i, \l}^{\frac{n-1}{n-2}}\varphi_{a_j, \l}^{\frac{n-1}{n-2}}dS_g
=
O_{A_{i, j}, \l}(\varepsilon_{i,j}^{\frac{n-1}{n-2}}\log \varepsilon_{i,j}),
\end{split}
\end{equation} 
where $A_{i, j}=(a_i, a_j)$, $\varepsilon_{i, j}:=\varepsilon_{i, j}(A, \bar\l)$ with $\bar \l:=(\l, \cdots, \l)$, and for the meaning of  $O_{A_{i, j}, \l}(1)$ and $\epsilon_{i, j}(A, \bar\l)$, see respectively Section \ref{eq:notpre} and \eqref{eq:varepsilonij}.
\end{lem}
\begin{pf}
Using \eqref{eq:standardbubbles}, \eqref{eq:varphialb}, \eqref{eq:decompval}-\eqref{eq:walest} and setting \;$\varphi_i=\varphi_{a_i, \l }$\; for \;$i=1, \cdots, p$, we have that for every $i=1, \cdots, p$, the following estimate holds
\begin{equation}\label{eq:betaint1}
\begin{split}
\varphi_{i} \leq C \left(\frac{\lambda}{1+\lambda^{2}d_{\hat g}^{2}(a_{i},\cdot)}\right)^{\frac{n-2}{2}}\;\;\text{on}\;\; \partial M
\end{split}\end{equation}
for some large positive constant independent of \;$a_i$\;and \;$\l$\; with \;$\hat g$\; denoting the Riemannian metric induced by \;$g$\; on \;$\partial M$. Hence, using \eqref{eq:betaint1}, we have for \;$c>0$\; and small that the following estimate holds
\begin{equation}\label{eq:betaint2}
\begin{split}
\oint_{\partial M} \varphi_{i}^{\frac{n-1}{n-2}} \varphi_{j}^{\frac{n-1}{n-2}}dS_g
\leq &
C \underset{B_{0}^{\R^{n-1}}(c)}{\int}\left(\frac{\lambda }{1+\lambda ^2|x|^{2}}\right)^{\frac{n-1}{2}}\left(\frac{ \lambda  }{1 +\lambda ^{2} 
d_{\hat g}^{2}(a_{j},\exp_{i}(x))}\right)^{\frac{n-1}{2}}dx\\
& +
C  \frac{1}{\lambda ^{\frac{n-1}{2}}}\underset{B_{0}^{\R^{n-1}}(c)}{\int}\left(\frac{ \lambda }{1+\lambda ^2|x|^{2}}\right)^{\frac{n-1}{2}}dx
+
O_{A_{i, j}, \l}\left(\frac{1}{\lambda ^{n-1}}\right) \\
= &
C \underset{B_{0}^{\R^{n-1}}(c\lambda)}{\int}\left(\frac{1}{1+r^{2}}\right)^{\frac{n-1}{2}}
\left(\frac{1}{1 
+
\lambda^{2}  d_{\hat g}^{2}(a_{j},\exp_{i}(\frac{x}{\lambda}))}\right)^{\frac{n-1}{2}}dx\\
& +
C  \frac{1}{\lambda^{n-1}} 
\underset{B_{0}^{\R^{n-1}}(c\l)}{\int}\left(\frac{1}{1+|x|^{2}}\right)^{\frac{n-1}{2}}dx
+
O_{A_{i, j}, \l}\left(\frac{1}{\lambda^{n-1}}\right),
\end{split}\end{equation} 
for some large positive constant \;$C$\; independent of \;$A_{i, j}$\; and \;$\l$\; with \;$r=|x|$\; and \;$\exp_i:=\hat{\exp}_{a_i}$ (for its meaning see Section \ref{eq:notpre}). Thus appealing to \eqref{eq:betaint2}, we infer that
\begin{equation}\label{eq:betaint3}
\begin{split}
\oint_{\partial M} \varphi_{i}^{\frac{n-1}{n-2}}  \varphi_{j}^{\frac{n-1}{n-2}}dS_g 
\leq &
C\underset{B_{0}^{\R^{n-1}}(c\l)}{\int}\left(\frac{1}{1+|x|^{2}}\right)^{\frac{n-1}{2}}\left(\frac{1}{1 +\lambda^{2}  
 d_{\hat g}^{2}(a_{j},\exp_{i}(\frac{x}{\lambda}))}\right)^{\frac{n-1}{2}}dx
+
O_{A_{i, j}, \l}(\frac{\log \lambda}{\lambda^{n-1}}).
\end{split}\end{equation} 
Thus \eqref{eq:betaint} follows from \eqref{eq:betaint3} if 
\begin{equation}\begin{split}
d_{\hat g}( a _{i}, a _{j})\geq 3c.
\end{split}\end{equation} 
Hence, to complete the proof of the lemma it remains to treat the case \;$d_{\hat g}( a _{i}, a _{j})<3c$. To do that, we set 
\begin{equation*}
\mathcal{B}=\{x\in \bar \R^{n-1}:\;\;\;\;\frac{1}{2}d_{\hat g}( a _{i}, a _{j})\leq \vert  \frac{x}{\lambda }\vert \leq 2 d_{\hat g}( a _{i}, a _{j})\},
\end{equation*}
and use \eqref{eq:betaint3} and the triangle inequality to get for \;$c>0$\; sufficiently small that the following estimate holds
\begin{equation}\begin{split}
\oint_{\partial M} \varphi_{i}^{\frac{n-1}{n-2}}&\varphi_{j}^{\frac{n-1}{n-2}}dS_g
\leq 
C \underset{\mathcal{B}}{\int}
\left(\frac{1}{1+|x|^{2}}\right)^{\frac{n-1}{2}}\left(\frac{1}{1 +\lambda^{2} d_{\hat g}^{2}( a _{j}, \exp_{i}( \frac{x}{\lambda }))}\right)^{\frac{n-1}{2}}dx +
O_{A_{i, j}, \l}(\varepsilon_{i,j}^{\frac{n-1}{n-2}}\log \varepsilon_{i,j})\\
\leq &
C
\left(\frac{1}{1+\vert\lambda d_{\hat g}( a _{i}, a _{j})\vert^{2}}\right)^{\frac{n-1}{2}}
\int_{\{\vert  \frac{x}{\lambda }\vert \leq 4 d_{\hat g}( a _{i}, a _{j})\}}
\left(\frac{1}{1+|x|^{2}}\right)^{\frac{n-1}{2}}dx 
+
O_{A_{i, j}, \l}(\varepsilon_{i,j}^{\frac{n-1}{n-2}}\log \varepsilon_{i,j})\\
= &
O_{A_{i, j}, \l}(\varepsilon_{i,j}^{\frac{n-1}{n-2}}\log \varepsilon_{i,j}),
\end{split}
\end{equation} 
where \;$C$\; is a large positive constant independent of \;$A_{i, j}$\; and \;$\l$, thereby completing the proof of the lemma.
\end{pf}

\section{Energy estimates for the barycenter technique}\label{eq:mappingbary}
In this section, we map \;$B_p(\partial M)$\; into some appropriate sublevels of the Euler-Lagrange functional \;$\mathcal{E}_g$\; via the Chen\cite{chen}'s bubbles. Precisely, we are going to derive sharp energy estimates for convex combinations of the bubbles \;$\varphi_{a, \l}$ given by \eqref{eq:varphialb} so that we can use them in the next section to run a suitable scheme of the  barycenter technique of Bahri-Coron\cite{bc}. In order to do that, we first make the following definition. For \;$p\in \N^*$, $\sigma:=\sum_{i=1}^p\alpha_i\d_{a_i}\in B_p(\partial M)$, and \;$\l\geq \frac{2}{\delta_0^{\frac{n-1}{2}}}$\; where $\delta_0$ is given by Lemma \ref{eq:brenchenenergy}, we define \;$f_p(\l): \;B_p(\partial M)\longrightarrow W^{1, 2}_+(\ov M)$ as follows 
\begin{equation}\label{eq:deffp}
f_p(\l)(\sigma):=\sum_{i=1}^p\alpha_i\varphi_{a_i, \l}.
\end{equation}
Now, we start the goal of this section with the following proposition which provides the first step to apply our scheme of the algebraic topological argument of Bahri-Coron\cite{bc}.
\begin{pro}\label{eq:baryest}
 There exists a large constant \;$C_0>0$, $\nu_0>1$\; and \;$0<\varepsilon_0\leq \delta_0$\; such that for every \;$p\in \N^*$ and every \;$0<\varepsilon\leq \varepsilon_0$, there exists \;$\l_p:=\l_p(\varepsilon):=\l_p(\nu_0, \varepsilon)\geq \frac{2}{\d_0^{\frac{n-1}{2}}}$\; such that for every \;$\l\geq \l_p$ and for every $\sigma=\sum_{i=1}^p\alpha_i\delta_{a_i}\in B_p(\partial M)$, we have\\\\
 1) If there exist \;$i_0\neq j_0$\; such that \;$\frac{\alpha_{i_0}}{\alpha_{j_0}}>\nu_0$\; or if \;$\sum_{i\neq j}\varepsilon_{i, j}> \varepsilon$, then
 $$
\mathcal{E}_g(f_p(\l))(\sigma))\leq p^{\frac{1}{n-1}}\mathcal{Q}(B^n),
$$ 
where \;$\mathcal{Q}(B^n)$\; is defined by \eqref{eq:yamabehsphere} and \;$\varepsilon_{i, j}:=\varepsilon_{i, j}(A, \bar \l)$ with $\bar \l:=(\l, \cdots, \l)$ and for the definition of \;$\varepsilon(A, \bar \l)$, see \eqref{eq:varepsilonij}.\\
2) If for every $i\neq j$ we have \;$\frac{\alpha_{i}}{\alpha_j}\leq\nu_0$\; and if \;$\sum_{i\neq j}\varepsilon_{i, j}\leq \varepsilon$, then
$$
\mathcal{E}_g(f_p(\l))(\sigma))\leq p^{\frac{1}{n-1}}\mathcal{Q}(B^n)\left(1+\frac{C_0}{\l^{n-2}}-c_g\frac{(p-1)}{\l^{n-2}}\right),
$$
where \;$c_g$\; is is defined by \eqref{eq:defcg}.
\end{pro}
\vspace{6pt}

\noindent
Like in \cite{martndia1}, Proposition \ref{eq:baryest} will be derived from  the following technical Lemma.
\begin{lem}\label{eq:baryestaux}
We have that the following holds:\\
1)
For every \;$\epsilon>0$\; and small and for every \;$p\in \N^*$, there exists\; $\lambda_{p}:=\l_p(\epsilon)\geq\frac{2}{\d_0^{\frac{n-1}{2}}}$ such that for every \;$\l\geq \l_p$ and for every $\sigma:=\sum_{i=1}^p\alpha_i\d_{a_i}\in B_p(\partial M)$, we have
\begin{equation*}
\begin{split}
\sum_{i\neq j}\epsilon_{i,j}>\epsilon
\end{split}
\end{equation*}
implies
\begin{equation*}
\begin{split}
\mathcal{E}_{g}(f_p(\l)(\sigma))< p^{\frac{1}{n-1}}\mathcal{Q}(B^n),
\end{split}
\end{equation*}
where \;$\epsilon_{i, j}:=\epsilon_{i, j}(A, \l)$\; is defined by \eqref{eq:intweight}.\\
2) For every \;$\nu>1$, for every \;$\epsilon>0$\;and small, and for every \;$p\in \N^*$, there exists \;$\lambda_{p}:=\lambda_{p}(\epsilon,\nu)\geq\frac{2}{\d_0^{\frac{n-1}{2}}}$\; such that for every \;$\l\geq \l_p$ and for every \;$\sigma:=\sum_{i}^p\alpha_i\d_{a_i}\in B_p(\partial M)$, we have
\begin{equation*}\begin{split}
\exists\;{i_0\neq j_0}\;\;\text{such that}\;\;\;\frac{\alpha_{i_0}}{\alpha_{j_0}} >\nu
\; \;\;\;\text{and}\; \;\;\;\sum_{i\neq j}\epsilon_{i,j}\leq\epsilon
\end{split}
\end{equation*}
imply
\begin{equation*}
\begin{split}
\mathcal{E}_{g}(f_p(\l)(\sigma))< p^{\frac{1}{n-1}}\mathcal{Q}(B^n).
\end{split}
\end{equation*}
3) There exists \;$C_0>0$, \;$\nu_0>1$, \;$\l_0\geq\frac{2}{\d_0^{\frac{n-1}{2}}}$\; and \;$0<\epsilon_0\leq\delta_0$\; such that for every \;$1<\nu\leq \nu_0$, for every \;$0<\epsilon\leq\epsilon_0$, for every \;$p\in \N^*$, for every \;$\l\geq \l_0$, and for every $\sigma:=\sum_{i=1}^p\alpha_i\d_{a_i}\in B_p(\partial M)$, we have
\begin{equation}\begin{split}
\frac{\alpha_{i}}{\alpha_{j}} \leq\nu\;\;\;\forall i, j,\;\;
\; \text{and}\;\;\; \sum_{i\neq j}\epsilon_{i,j}\leq\epsilon
\end{split}\end{equation}
imply
\begin{equation}\begin{split} 
\mathcal{E}_{g}(f_p(\l)(\s))
\leq &
p^{\frac{1}{n-1}}\mathcal{Q}(B°)
\left(
1+\frac{C_0}{\lambda^{n-2}}
-
c_g\frac{(p-1)}{\lambda^{n-2}}
\right).
\end{split}\end{equation} 
\end{lem}
\begin{pf}
The strategy of the proof is the same as the one of Lemma 3.2 in \cite{martndia1}. For the sake of completeness we will provide full details. First of all, we set 
\begin{equation}\label{eq:numden}
 \mathcal{N}_g(u):=\langle L_gu,u\rangle+\langle B_gu, u\rangle, \;\;\mathcal{D}_g(u):=\left(\oint_{\partial M} u^{\frac{2(n-1)}{n-2}}dS_g\right)^{\frac{n-2}{n-1}}, \;\;\;\;u\in W^{1, 2}_+(\ov M),
\end{equation}
and use \eqref{eq:escobarfunctional} to have
\begin{equation}\label{eq:eulerlagrange}
\mathcal{E}_g(u)=\frac{\mathcal{N}_g(u)}{\mathcal{D}_g(u)},\;\;\;\;u\in W^{1, 2}_+(\ov M).
\end{equation}
Furthermore, for $p\in \N^*$, $\sigma:=\sum_{i=1}^p\alpha_i\d_{a_i}\in B_p(\partial M)$\; and \;$\l\geq\frac{2}{\d_0^{\frac{n-1}{2}}}$, we set (as in the proof of Lemma \ref{eq:interactestl})
\begin{equation}\label{eq:varphii1}
 \varphi_i=\varphi_{a_i, \l}, \;\;i=1, \cdots, p.
\end{equation}
Now, we start with the proof of point 1). To do so, we first use Lemma \ref{eq:interactestl}, \eqref{eq:deffp}, \eqref{eq:numden}, \eqref{eq:varphii1}, and H\"older's inequality to estimate\; $\mathcal{N}_g(f_p(\l)(\s))$\; as follows
\begin{equation}\label{eq:ngest1}
\begin{split} 
\mathcal{N}_{g}(f_p(\l)(\s))
=&
c_{0}(1+o_{A, \bar\alpha, \l}(1))
\oint_{\partial M} \left(\sum_{i=1}^p\alpha_{i}\varphi_{i}^{\frac{n}{n-2}}\right)\left(\sum_{j=1}^p\alpha_{j}\varphi_{j}\right)dS_g\\
= &
c_{0}(1+o_{A, \bar\alpha, \l}(1))
\oint_{\partial M} \left(\frac{\sum_{i=1}^p\alpha_{i}\varphi_{i}^{\frac{n}{n-2}}}{ \sum_{j=1}^p\alpha_{j}\varphi_{j}}\right)\left(\sum_{j=1}^p\alpha_{j}\varphi_{j}\right)^{2}dS_g \\
\leq &
c_{0}(1+o_{A, \bar\alpha, \l}(1))
D_{g}(u)\Vert\frac{\sum_{i=1}^p\alpha_{i}\varphi_{i}^{\frac{n}{n-2}}}{\sum_{j=1}^p\alpha_{j}\varphi_{j}}\Vert_{L^{n-1}(\partial M)},
\end{split}
\end{equation} 
where\; $A:=(a_1, \cdots, a_p)$, $\bar\alpha:=(\alpha_1, \cdots, \alpha_p)$ and for the meaning of \;$o_{A, \bar\alpha, \l}(1)$, see Section \ref{eq:notpre}.
Thus, using the convexity of the map $x\longrightarrow x^{\beta}$ with $\beta>1$, we derive that \eqref{eq:ngest1} implies  
\begin{equation}\label{eq:ngest2}
\begin{split} 
\mathcal{N}_{g}(f_p(\l)(\s))
\leq &
c_{0}(1+o_{A, \bar\alpha, \l}(1))
\mathcal{D}_{g}(u)\left(\oint_M\left(\sum_{i=1}^p\frac{\alpha_{i}\varphi_{i}}{\sum_{j=1}^p\alpha_{j}\varphi_{j}}\varphi_{i}^{\frac{2}{n-2}}\right)^{n-1}\right)^{\frac{1}{n-1}}dS_g \\
\leq &
c_{0}(1+o_{A, \bar\alpha, \l}(1))
\mathcal{D}_{g}(u)\left(\sum_{i=1}^p\oint_{\partial M}\frac{\alpha_{i}\varphi_{i}}{\sum_{j=1}^p\alpha_{j}\varphi_{j}}\varphi_{i}^{\frac{2(n-1)}{n-2}}\right)^{\frac{1}{n-1}}dS_g.
\end{split}\end{equation} 
Hence, clearly Lemma \ref{eq:bsvolestl}, \eqref{eq:eulerlagrange} and \eqref{eq:ngest2} imply for any pair \;$i\neq j$ ($i, j=1, \cdots, p$)
\begin{equation}\label{eq:ngest3}
\begin{split} 
\mathcal{E}_{g}(f_p(\l)(\s))
\leq & 
c_{0}(1+o_{A, \bar\alpha, \l}(1))
\left(
c_{1}(p-1)
+
\oint_{\partial M} \frac{\alpha_{i}\varphi_{i}}{\alpha_{i}\varphi_{i}+\alpha_{j}\varphi_{j}}\varphi_{i}^{\frac{2(n-1)}{n-2}}dS_g
\right)^{\frac{1}{n-1}} \\
\leq &
c_{0}(1+o_{A, \bar\alpha,  \l}(1))
\left(
c_{1}p
-
\oint_{\partial M} \frac{\alpha_{j}\varphi_{j}}{\alpha_{i}\varphi_{i}+\alpha_{j}\varphi_{j}}\varphi_{i}^{\frac{2(n-1)}{n-2}}
dS_g\right)^{\frac{1}{n-1}}, \\
\end{split}\end{equation} 
and we may assume \;$\alpha_{i}\leq\alpha_{j}$\; by symmetry. Now, we are going to estimate from below the quantity $\oint_{\partial M} \frac{\alpha_{j}\varphi_{j}}{\alpha_{i}\varphi_{i}+\alpha_{j}\varphi_{j}}\varphi_{i}^{\frac{2(n-1)}{n-2}}dS_g
$. In order to do that, for \;$\gamma>0$, we set
\begin{equation}\label{eq:defmathaij}
\begin{split} 
\mathcal{A}_{i,j}=
\{x\in \partial M:\;\;
\varphi_{i}(x)\geq \gamma (\frac{\alpha_{i}}{\alpha_{j}}\varphi_{i}(x)+\varphi_{j}(x))
\},
\end{split}\end{equation} 
and use \eqref{eq:defmathaij} to have
\begin{equation}\label{eq:ngest4}
\begin{split} 
\oint_{\partial M} \frac{\alpha_{j}\varphi_{j}\varphi_{i}^{\frac{2(n-1)}{n-2}}}{\alpha_{i}\varphi_{i}+\alpha_{j}\varphi_{j}}dS_g
\geq &
\oint_{\mathcal{A}_{i,j}} \frac{\varphi_{j}}{\frac{\alpha_{i}}{\alpha_{j}}\varphi_{i}+\varphi_{j}}\varphi_{i}^{\frac{2(n-1)}{n-2}}dS_g
\geq 
\gamma 
\oint_{\mathcal{A}_{i,j}}
\varphi_{i}^{\frac{n}{n-2}}\varphi_{j}dS_g \\
= &
\gamma 
\left(
\oint_{\partial M} \varphi_{i}^{\frac{n}{n-2}}\varphi_{j}dS_g-\oint_{\mathcal{A}_{i,j}^{c}}\varphi_{i}^{\frac{n}{n-2}}\varphi_{j}dS_g
\right) \\
\geq &
\gamma 
\left(
\oint_{\partial M} \varphi_{i}^{\frac{n}{n-2}}\varphi_{j}dS_g
-
\gamma ^{\frac{2}{n-2}}
\oint_{\mathcal{A}_{i,j}^{c}}\left(\frac{\alpha_{i}}{\alpha_{j}}\varphi_{i}+\varphi_{j}\right)^{\frac{2}{n-2}}\varphi_{i}\varphi_{j}dS_g
\right),
\end{split}\end{equation} 
where \;$\mathcal{A}_{i,j}^c:=\partial M\setminus \mathcal{A}_{i, j}$.
Next, since \;$\frac{\alpha_{i}}{\alpha_{j}}\leq 1$, then appealing to \eqref{eq:ngest4}, we infer that the following estimate holds 
\begin{equation}\label{eq:ngest5}
\begin{split} 
\oint_{\partial M} \frac{\alpha_{j}\varphi_{j}\varphi_{i}^{\frac{2(n-1)}{n-2}}}{\alpha_{i}\varphi_{i}+\alpha_{j}\varphi_{j}}dS_g
\geq &
\gamma 
\left(
\oint_{\partial M} \varphi_{i}^{\frac{n}{n-2}}\varphi_{j}dS_g
-
C\gamma ^{\frac{2}{n-2}}
\oint_{\partial M}(\varphi_{i}^{\frac{2}{n-2}}+\varphi_{j}^{\frac{2}{n-2}})\varphi_{i}\varphi_{j}dS_g
\right),
\end{split}\end{equation} 
for some large positive constant \;$C$\; independent of \;$A$, $\l$ and $\gamma$. Thus, ii) of point 1) of Lemma \ref{eq:interactestl} and \eqref{eq:ngest5} imply that for $\gamma>0$ sufficiently small, there holds
\begin{equation}\label{eq:ngest6}
\begin{split} 
\oint_{\partial M} \frac{\alpha_{j}\varphi_{j}}{\alpha_{i}\varphi_{i}+\alpha_{j}\varphi_{j}}\varphi_{i}^{\frac{2(n-1)}{n-2}}dS_g
\geq &
\frac{\gamma }{2}
\int_{\partial M} \varphi_{i}^{\frac{n}{n-2}}\varphi_{j}dS_g.
\end{split}\end{equation} 
Hence, combining \eqref{eq:ngest3} and \eqref{eq:ngest6}, we conclude that for any pair \;$i\neq j$, the following estimate holds
\begin{equation}\label{eq:ngest7}
\begin{split} 
\mathcal{E}_{g}(f_p(\l)(\s))
\leq & 
\left(1+o_{A, \bar\alpha, \l}(1)\right)\mathcal{Q}(B^n)
\left(
p
-
\frac{\gamma }{2c_1}
\oint_{\partial M} \varphi_{i}^{\frac{n}{n-2}}\varphi_{j}dS_g
\right)^{\frac{1}{n-1}}.
\end{split}\end{equation} 
Clearly \eqref{eq:ngest7} implies, that we always have 
\begin{equation}\begin{split} 
\mathcal{E}_g(f_p(\l)(\s))\leq \left(1+o_{A, \bar\alpha, \l}(1)\right)p^{\frac{1}{n-1}}\mathcal{Q}(B^n)
\end{split}\end{equation} 
and in case $\sum_{i\neq j}\epsilon_{i,j}> \epsilon$ 
\begin{equation}\begin{split} 
\mathcal{E}_{g}(f_p(\l)(\s))
\leq & 
\left(1+o_{A, \bar\alpha, \l}(1)\right)p^{\frac{1}{n-1}}\mathcal{Q}(B^n)
\left(
1
-
\frac{\gamma \epsilon}{2pc_1}
\right)^{\frac{1}{n-1}}.
\end{split}\end{equation} 
thereby ending the proof of point 1). Now, we are going to treat the second case. Hence, we may assume
\begin{equation}\begin{split} 
\sum_{i\neq j}\epsilon_{i,j}\ll 1
\end{split}\end{equation} 
and thus according to Lemma \ref{eq:interactestl}
\begin{equation}\label{eq:varepijep}
\begin{split}
\epsilon_{i,j}=(1+o_{\varepsilon_{i, j}}(1))\varepsilon_{i,j}
\;\text{ and }\;
\lambda d_{\hat g}(a_{i},a_{j})\gg 1,
\end{split}
\end{equation}
and for the meaning of \;$o_{\epsilon_{i, j}}(1)$, see Section \ref{eq:notpre}. We then use Lemma \ref{eq:interactestl},  \eqref{eq:numden}, and \eqref{eq:varepijep} to have 
\begin{equation}\begin{split} \label{eq:ngest8}
\mathcal{N}_{g}(f_p(\l)(\s))
= &
\sum_{i=1}^p\sum_{j=1}^p\alpha_{i}\alpha_{j}\left(\langle L_{g}\varphi_{i}, \varphi_{j}\rangle+\langle B_g \varphi_i, \varphi_j\rangle\right) \\
= &
\sum_{i=1}^p\alpha_{i}^{2}\left(\langle L_{g}\varphi_{i}, \varphi_{i}\rangle+\langle B_g \varphi_i, \varphi_i \rangle \right)
+
\sum_{i\neq j}\alpha_{i}\alpha_{j}\left(\langle L_{g}\varphi_{i}, \varphi_{j}\rangle+\langle B_g\varphi_i, \varphi_j\rangle \right) \\
= &
\sum_{i}\alpha_{i}^{2}\mathcal{E}_{g}(\varphi_{i})\mathcal{D}_{g}(\varphi_{i})
+
c_{0}(1+o_{\sum_{i,\neq j}\varepsilon_{i, j}}(1))\sum_{i\neq j}\alpha_{i}\alpha_{j}\varepsilon_{i,j}
\end{split}\end{equation} 
and (for the meaning of \;$o_{\sum_{i\neq j}\epsilon_{i, j}}(1)$, see Section \ref{eq:notpre})
\begin{equation}\label{eq:dgest1}
\begin{split} 
\mathcal{D}^{\frac{n-1}{n-2}}_{g}(f_p(\l)(\s))
= &
\oint_{\partial M} \left(\sum_{i=1}^p\alpha_{i}\varphi_{i}\right)^{\frac{2(n-1)}{n-2}}dS_g
=
\sum_{i=1}^p\alpha_{i}\oint_{\partial M} \left(\sum_{j=1}^p\alpha_{j}\varphi_{j}\right)^{\frac{n}{n-2}}\varphi_{i}dS_g \\
= &
\sum_{i=1}^p\alpha_{i}\oint_{\partial M} \left(\alpha_{i}\varphi_{i}+\sum_{j=1, \;j\neq i}^{p}\alpha_{j}\varphi_{j}\right)^{\frac{n}{n-2}}\varphi_{i}dS_g.
\end{split}\end{equation} 
To proceed further, we set 
$
\mathcal{A}_{i}=\{x\in \partial M:\;\;\alpha_{i}\varphi_{i}(x)>\sum_{j=1, \;j\neq i}^{p}\alpha_{j}\varphi_{j}(x)\},
$
and use Taylor expansion to obtain
\begin{equation}\label{eq:dgest2}
\begin{split} 
\mathcal{D}^{\frac{n-1}{n-2}}_{g}(f_p(\l)(\s))
= &
\sum_{i=1}^p\alpha_{i}^{\frac{2(n-1)}{n-2}}\oint_{\mathcal{A}_{i}} \varphi_{i}^{\frac{2(n-1)}{n-2}}dS_g
+
\frac{n}{n-2}\sum_{i\neq j}\alpha_{i}^{\frac{n}{n-2}}\alpha_{j}\oint_{\mathcal{A}_{i}} \varphi_{i}^{\frac{n}{n-2}}\varphi_{j}dS_g\\
+ &
\sum_{i=1}^p\alpha_{i}\oint_{\mathcal{A}_{i}^{c}}\left(\sum_{ j=1, \;j\neq i}^{p}\alpha_{j}\varphi_{j}\right)^{\frac{n}{n-2}}\varphi_{i}dS_g +
O_{A, \bar\alpha, \l}
\left(\sum_{i\neq j}\alpha_{i}^{\frac{n-1}{n-2}}\alpha_{j}^{\frac{n-1}{n-2}}\oint_{\partial M} \varphi_{i}^{\frac{n-1}{n-2}}\varphi_{j}^{\frac{n-1}{n-2}}dS_g\right),
\end{split}\end{equation} 
where \;$\mathcal{A}_i^c:=\partial M\setminus\mathcal{A}_i$, $A:=(a_1, \cdots, a_p)$, $\bar\alpha:=(\alpha_1, \cdots, \alpha_p)$, $O_{A, \bar\alpha, \l}(1)$ is defined as in Section \ref{eq:notpre}, and we made use of \;$n\geq 3$\; and the algebraic relation
\begin{equation}\begin{split}
(a+b)^{\frac{n-1}{n-2}}\leq C_n(a^{\frac{n-1}{n-2}}+b^{\frac{n-1}{n-2}})
\end{split}\end{equation}
for \;$a,b\geq 0$\; and \;$C_n$\; a positive constant depending only on \;$n$. Moreover, since
\begin{equation}\begin{split} 
(a+b)^{\frac{n}{n-2}}\geq a^{\frac{n}{n-2}}+b^{\frac{n}{n-2}},
\end{split}\end{equation} 
then \eqref{eq:dgest2} implies
\begin{equation}\label{eq:dgest3}
\begin{split} 
\mathcal{D}^{\frac{n-1}{n-2}}_{g}(f_p(\l)(\s))
= &
\sum_{i=1}^p\alpha_{i}^{\frac{2(n-1)}{n-2}}\oint_{\partial M} \varphi_{i}^{\frac{2(n-1)}{n-2}}dS_g 
+
\frac{n}{n-2}\sum_{i\neq j}\alpha_{i}^{\frac{n}{n-2}}\alpha_{j}\oint_{\partial M} \varphi_{i}^{\frac{n}{n-2}}\varphi_{j}dS_g \\
& + 
\sum_{i\neq j}\alpha_{i}\alpha_{j}^{\frac{n}{n-2}}\oint_{A_{i}^{c}} \varphi_{j}^{\frac{n}{n-2}}\varphi_{i}dS_g \\
& +
O_{A, \l}
\left(\sum_{i\neq j}\alpha_{i}^{\frac{n-1}{n-2}}\alpha_{j}^{\frac{n-1}{n-2}}\oint_{\partial M} \varphi_{i}^{\frac{n-1}{n-2}}\varphi_{j}^{\frac{n-1}{n-2}}dS_g\right) \\
= &
\sum_{i=1}^p\alpha_{i}^{\frac{2(n-1)}{n-2}}\oint_{\partial M} \varphi_{i}^{\frac{2(n-1)}{n-2}}dS_g 
+
2\frac{n-1}{n-2}\sum_{i\neq j}\alpha_{i}^{\frac{n}{n-2}}\alpha_{j} \oint_{\partial M} \varphi_{i}^{\frac{n}{n-2}}\varphi_{j}dS_g
 \\
& +
O_{A, \bar\alpha, \l}
\left(\sum_{i\neq j}\alpha_{i}^{\frac{n-1}{n-2}}\alpha_{j}^{\frac{n-1}{n-2}}\oint_{\partial M} \varphi_{i}^{\frac{n-1}{n-2}}\varphi_{j}^{\frac{n-1}{n-2}}dS_g\right).
\end{split}\end{equation} 
So, using Lemma \ref{eq:alphabetainteract} and \eqref{eq:varepijep}, we have that \eqref{eq:dgest3} implies
\begin{equation}\label{eq:dgest4}
\begin{split} 
\mathcal{D}^{\frac{n-1}{n-2}}_{g}(f_p(\l)(\s))
= &
\sum_{i=1}^p\alpha_{i}^{\frac{2(n-1)}{n-2}}\oint_{\partial M} \varphi_{i}^{\frac{2(n-1)}{n-2}}dS_g  \\
& +
\frac{2(n-1)}{n-2}\left(1+o_{\sum_{i\neq j}\varepsilon_{i, j}}(1)\right)\sum_{i\neq j}\alpha_{i}^{\frac{n}{n-2}}\alpha_{j} \varepsilon_{i,j}
+o_{\sum_{i\neq j}\varepsilon_{i, j}}\left(\sum_{i\neq j}\alpha_{i}^{\frac{n-1}{n-2}}\alpha_{j}^{\frac{n-1}{n-2}}\varepsilon_{i,j}\right),
\end{split}
\end{equation} 
Thus, using Young's inequality and the symmetry of $\varepsilon_{i,j}$, we infer from \eqref{eq:dgest4} that the following estimate holds
\begin{equation}\label{eq:dgest5}
\begin{split} 
\mathcal{D}^{\frac{n-1}{n-2}}_{g}(f_p(\l)(\s)))
= &
\sum_{i=1}^p\alpha_{i}^{\frac{2(n-1)}{n-2}}\oint_{\partial M} \varphi_{i}^{\frac{2(n-1)}{n-2}} dS_g 
+
\frac{2(n-1)}{n-2}(1+o_{\sum_{i\neq j}\varepsilon_{i, j}}(1))\sum_{i\neq j}\alpha_{i}^{\frac{n}{n-2}}\alpha_{j} \varepsilon_{i,j}.
\end{split}\end{equation} 
Hence, using again Young's inequality, Taylor expansion, and Lemma \ref{eq:bsvolestl}, we have that \eqref{eq:dgest5} gives
\begin{equation}\label{eq:dgest6}
\begin{split} 
\mathcal{D}_{g}(f_p(\l)(\s))
=& 
\left(\sum_{i=1}^p\alpha_{i}^{\frac{2(n-1)}{n-2}}\oint_{\partial M} \varphi_{i}^{\frac{2(n-1)}{n-2}}dS_g\right)^{\frac{n-2}{n-1}}\times
\left(
1
+
\frac
{
\frac{2(n-1)}{n-2}\left(1+o_{\sum_{i\neq j}\varepsilon_{i, j}}(1)\right)\sum_{i\neq j}\alpha_{i}^{\frac{n}{n-2}}\alpha_{j} \varepsilon_{i,j}
}
{
\sum_{i=1}^p\alpha_{i}^{\frac{2(n-1)}{n-2}}\oint_{\partial M} \varphi_{i}^{\frac{2(n-1)}{n-2}}dS_g
}
\right)^{\frac{n-2}{n-1}} 
 \\
= &
\left(\sum_{i=1}^p\alpha_{i}^{\frac{2(n-1)}{n-2}}\oint_{\partial M} \varphi_{i}^{\frac{2(n-1)}{n-2}}dS_g\right)^{\frac{n-2}{n-1}}
+
2
\frac
{
\left(1+o_{\sum_{i\neq j}\varepsilon_{i, j}}(1)\right)\sum_{i\neq j}\alpha_{i}^{\frac{n}{n-2}}\alpha_{j} \varepsilon_{i,j}
}
{
\left(\sum_{i=1}^p\alpha_{i}^{\frac{2(n-1)}{n-2}}\oint_{\partial M} \varphi_{i}^{\frac{2(n-1)}{n-2}}dS_g\right)^{\frac{1}{n-1}}
}
\\
= &
\left(\sum_{i=1}^p\alpha_{i}^{\frac{2(n-1)}{n-2}}\mathcal{D}_{g}^{\frac{n-1}{n-2}}(\varphi_{i})\right)^{\frac{n-2}{n-1}}
+
2c_{1}^{-\frac{1}{n-1}}
\frac
{
\left(1+o_{\sum_{i\neq j}\varepsilon_{i, j}}(1)\right)\sum_{i\neq j}\alpha_{i}^{\frac{n}{n-2}}\alpha_{j} \varepsilon_{i,j}
}
{
\left(\sum_{i=1}^p\alpha_{i}^{\frac{2(n-1)}{n-2}}\right)^{\frac{1}{n-1}}
}.
\end{split}\end{equation} 
Now, combining \eqref{eq:ngest8} and  \eqref{eq:dgest6}, and using again Taylor expansion,  we obtain 
\begin{equation}\label{eq:ngest9}
\begin{split} 
\mathcal{E}_{g}(f_p(\l)(\s))
\leq &
\frac
{
\sum_{i=1}^p\alpha_{i}^{2}\mathcal{E}_{g}(\varphi_{i})\mathcal{D}_{g}(\varphi_{i})
+
c_{0}\left(1+o_{\sum_{i\neq j}\varepsilon_{i, j}}(1)\right)\sum_{i\neq j}\alpha_{i}\alpha_{j}\varepsilon_{i,j}
}
{
\left(\sum_{i=1}^p\alpha_{i}^{\frac{2(n-1)}{n-2}}\mathcal{D}_{g}^{\frac{n-1}{n-2}}(\varphi_{i})\right)^{\frac{n-2}{n-1}}
+
2c_{1}^{-\frac{1}{n-1}}
\frac
{
\left(1+o_{\sum_{i\neq j}\varepsilon_{i, j}}(1)\right)\sum_{i\neq j}\alpha_{i}^{\frac{n}{n-2}}\alpha_{j} \varepsilon_{i,j}
}
{
\left(\sum_{i=1}^p\alpha_{i}^{\frac{2(n-1)}{n-2}}\right)^{\frac{1}{n-1}}
}
} \\
= &
\frac
{
\sum_{i=1}^p\mathcal{E}_{g}(\varphi_{i})\alpha_{i}^{2}\mathcal{D}_{g}(\varphi_{i})
}
{
\left(\sum_{i=1}^p\alpha_{i}^{\frac{2(n-1)}{n-2}}\mathcal{D}_{g}^{\frac{n-1}{n-2}}(\varphi_{i})\right)^{\frac{n-2}{n-1}}
} +
\frac{c_{0}}{c_{1}^{\frac{n-2}{n-1}}}\frac{\left(1+o_{\sum_{i\neq j}\varepsilon_{i, j}}(1)\right)\sum_{i\neq j}\alpha_{i}\alpha_{j}\varepsilon_{i,j}}{\left(\sum_{i=1}^p\alpha_{i}^{\frac{2(n-1)}{n-2}}\right)^{\frac{n-2}{n-1}}} \\
& -
\frac{2c_{0}}{c_{1}^{\frac{n-2}{n-1}}}\frac{\left(1+o_{\sum_{i\neq j}\varepsilon_{i, j}}(1)\right)\left(\sum_{i=1}^p\alpha_{i}^{2}\right)\left(\sum_{i\neq j}\alpha_{i}^{\frac{n}{n-2}}\alpha_{j}\varepsilon_{i,j}\right)}{(\sum_{i=1}^p\alpha_{i}^{\frac{2(n-1)}{n-2}})^{\frac{2n-3}{n-1}}}.
\end{split}\end{equation} 
Hence, using \eqref{eq:relationcy} and rearranging the terms in \eqref{eq:ngest9}, we get
\begin{equation}\label{eq:ngest10}
\begin{split}
\mathcal{E}_{g}(f_p(\l)(\s))
\leq &
\max_{i=1, \cdots, p}\mathcal{E}_{g}(\varphi_{i})\frac
{
\sum_{i=1}^p\alpha_{i}^{2}\mathcal{D}_{g}(\varphi_{i})
}
{
\left(\sum_{i=1}^p\left(\alpha_{i}^{2}\mathcal{D}_{g}(\varphi_{i})\right)^{\frac{n-1}{n-2}}\right)^{\frac{n-2}{n-1}}
}
\\ 
& +
\frac
{\left(1+o_{\sum_{i\neq j}\varepsilon_{i, j}}(1)\right)\mathcal {Q}(B^n)}
{\left(\sum_{i=1}^p\alpha_{i}^{\frac{2(n-1)}{n-2}}\right)^{\frac{n-2}{n-1}}}\left(
\sum_{i\neq j}\left[1-2\frac{\alpha_{i}^{\frac{2}{n-2}}\left(\sum_{i=1}^p\alpha_{i}^{2}\right)}{\left(\sum_{i=1}^p\alpha_{i}^{\frac{2(n-1)}{n-2}}\right)}\right]\alpha_{i}\alpha_{j}\frac{\varepsilon_{i,j}}{c_1}\right).
\end{split}\end{equation} 
This inequality has the following impact. First note, that the function
\begin{equation}\label{eq:defgammafunct}
\begin{split} 
\Gamma:\{\gamma \in [0,1]^{p}:&\;\;\;  \sum_{i=1}^p\gamma_{i} =1\}\longrightarrow  \mathbb{R}_{+}\\&\gamma
\longrightarrow 
\frac{\sum_{i=1}^p\gamma_{i}^{2}}{(\sum_{i=1}^p\gamma_{i}^{\frac{2(n-1)}{n-2}})^{\frac{n-2}{n-1}}}
\end{split}\end{equation} 
has the strict global maximum 
\begin{equation}\label{eq:gammamax}
\gamma_{\max}=\left(\frac{1}{p},\ldots,\frac{1}{p}\right)
\end{equation} 
with \;$\Gamma(\gamma_{\max})=p^{\frac{1}{n-1}}$. Thus, using Lemma \ref{eq:bsvolestl}, Lemma \ref{eq:bubbleestl}, \eqref{eq:numden}, \eqref{eq:varepijep}, and \eqref{eq:ngest10}, we infer that for any \;$\nu>0$, for every \;$\epsilon>0$\;and small, and for every \;$p\in\N^*$, there exists \;$\l_p:=\l_p(\nu, \epsilon)\geq\frac{2}{\d_0^{\frac{n-1}{2}}}$\; such for every \;$\l\geq \l_p$\; and for every $\sigma:=\sum_{i=1}^p\alpha_i\d_{a_i}\in B_p(\partial M)$, there holds
\begin{equation}\begin{split} 
\mathcal{E}_{g}(f_p(\l)(\s))<p^{\frac{1}{n-1}}\mathcal{Q}(B^n),
\end{split}\end{equation} 
whenever 
\begin{equation}\begin{split} 
\exists i_o\neq j_0 \;\;\;\text{such that}\;\;\; \frac{\alpha_{i_0}}{\alpha_{j_0}}>\nu\;
\;\;\;\; \text{and}\;\;\;\; 
\sum_{i\neq j}\varepsilon_{i,j}\leq\epsilon,
\end{split}\end{equation} 
thereby ending the proof of point 2). Now, we are going to treat point 3) and end the proof of the Lemma. Thus, we may assume
\begin{equation}\label{eq:assump1}
\begin{split} 
\forall\;\, i,j\;\;\frac{\alpha_{i}}{\alpha_{j}}= 1+o^+_{\mu}(1)\;\;\;\;
\;\text{and}\;\;\; \sum_{i\neq j}\epsilon_{i,j}\ll 1,
\end{split}\end{equation} 
where \;$o_{\mu}^+(1)$\; is a positive quantity depending only \;$\mu$\; with \;$\mu$\; small and verifying the property that it tends to \;$0$\; as \;$\mu$\; tends to \;$0$. So, using \eqref{eq:ngest10}, \eqref{eq:assump1}, and the properties of \;$\Gamma$\; (see \eqref{eq:defgammafunct} and  \eqref{eq:gammamax}), we infer that the following estimate holds
\begin{equation}\label{eq:ngest11}
\begin{split} 
\mathcal{E}_{g}(f_p(\l)(\s))
\leq &
\max_{i=1, \cdots, p}\mathcal{E}_{g}(\varphi_{i})p^{\frac{1}{n-1}}\\&
-(1+o_{\sum_{i\neq j}\varepsilon_{i, j}}(1)+o_{\mu}(1))
\frac
{\mathcal{Q}(B^n)}
{(\sum_{i=1}^p\alpha_{i}^{\frac{2(n-1)}{n-2}})^{\frac{n-2}{n-1}}}
\sum_{i\neq j}\alpha_{i}\alpha_{j}\frac{\varepsilon_{i,j}}{c_1} \\
= &
\max_{i=1, \cdots, p}\mathcal{E}_{g}(\varphi_{i})p^{\frac{1}{n-1}}
-(1+o_{\sum_{i\neq j}\varepsilon_{i, j}}(1)+o_{\mu}(1))
\mathcal{Q}(B^n)p^{\frac{2-n}{n-1}}
\sum_{i\neq j}\frac{\varepsilon_{i,j}}{c_1}.
\end{split}\end{equation} 
Now, using Lemma \ref{eq:bubbleestl}, \eqref{eq:varepijep},  \eqref{eq:assump1}, and \eqref{eq:ngest11}, we have that there exists \;$C_0>0$, $\nu_0>1$, $\l_0\geq \frac{2}{\delta_0^{\frac{n-1}{2}}}$\; and \;$0<\epsilon_0\leq\delta_0$\; such that for every \;$1<\nu\leq \nu_0$, for every \;$0<\epsilon\leq\epsilon_0$, for every \;$p\in \N^*$, for every \;$\l\geq \l_0$, and for every \;$\sigma:=\sum_{i=1}^p\alpha_i\d_{a_i}\in B_p(\partial M)$, we have if $\frac{\alpha_{i}}{\alpha_{j}} \leq\nu$\;\;\;$\forall i, j$\;\;
\; and\;\;\; $\sum_{i\neq j}\epsilon_{i,j}\leq\epsilon$, then there holds
\begin{equation}\label{eq:nges12}
\begin{split} 
\mathcal{E}_{g}(f_p(\l)(\s))
\leq &
p^{\frac{1}{n-1}}\mathcal{Q}(
B^n)
\left(
1+\frac{C_0}{\lambda^{n-2}}
-
\frac{1}{2c_1p}
\sum_{i\neq j}\varepsilon_{i,j}
\right).
\end{split}\end{equation} 
Thus, recalling that (see \eqref{eq:48}) 
\begin{equation}\begin{split}
\varepsilon_{i,j} = (1+o_{\varepsilon_{i, j}}(1))c_{3}\frac{G(a_{i}, a_{j})}{\lambda^{n-2}},
\end{split}\end{equation}
and using again \eqref{eq:varepijep}, we infer from \eqref{eq:nges12} that up to taking $\epsilon_0$ smaller, for every $1<\nu\leq \nu_0$, for every $0<\epsilon\leq\epsilon_0$, for every \;$p\in \N^*$, for every $\l\geq \l_0$, and for every $\sigma:=\sum_{i=1}^p\alpha_i\d_{a_i}\in B_p(\partial M)$, there holds  $$\frac{\alpha_{i}}{\alpha_{j}} \leq\nu\;\;\;\forall i, j\;\;
\; \text{and}\;\;\; \sum_{i\neq j}\epsilon_{i,j}\leq\epsilon$$ imply 
\begin{equation}\begin{split} 
\mathcal{E}_{g}(f_p(\l)(\s))
\leq &
p^{\frac{1}{n-1}}\mathcal{Q}(B^n)\left
(
1+\frac{C_0}{\lambda^{n-2}}
-
\frac{c_{3}}{4c_1p\lambda^{n-2}}
\sum_{i\neq j}G(a_{i},a_{j})
\right)\\
\leq &
p^{\frac{1}{n-1}}\mathcal{Q}(B^n)\left(
1+\frac{C_0}{\lambda^{n-2}}
-
c_g\frac{(p-1)}{\lambda^{n-2}}
\right),
\end{split}\end{equation} 
thereby ending the proof of point 3), and hence of the Lemma.
\end{pf}
\vspace{4pt}

\noindent
\begin{pfn}{ of Proposition \ref{eq:baryest}}\\
It follows from Lemma \ref{eq:baryestaux} by taking \;$C_0$\; and \;$\nu_0$\; to be the ones given by Lemma \ref{eq:baryestaux}, while \;$\varepsilon_0:=\frac{\epsilon_0}{2}$\;, and \;$\l_p:=\l_p(\varepsilon, \nu_0):=\max\{\l_p(\frac{\varepsilon}{2}), \l_p( 2\varepsilon, \nu_0), \l_0\}$, where\; $\epsilon_0$, $\l_p(\frac{\varepsilon}{2})$, $\l_p(2\varepsilon, \nu_0)$, and \;$\l_0$\; are as in Lemma \ref{eq:baryestaux}.
\end{pfn}
\vspace{6pt}

\noindent
Now, using Proposition \ref{eq:baryest}, we have the following corollary  which will be used together with Proposition \ref{eq:baryest} in the next section to carry a suitable algebraic topological argument of Bahri-Coron\cite{bc}.
\begin{cor}\label{eq:largep}
There exists $p_0\in \N^*$ large enough such that for every \;$0<\varepsilon\leq\varepsilon_0$, and for every for every \;$\l\geq\l_{p_0}$ (where \;;$\varepsilon_0$\; and \;$\l_{p_0}$\; are given by Proposition \ref{eq:baryest}), there holds 
$$
\mathcal{E}_g(f_{p_0}(\l)(B_{p_0}(\partial M)))\subset W_{p_0-1}.
$$
\end{cor}
\begin{pf}
It follows directly from Proposition \ref{eq:baryest} and the definition of \;$W_{p_0-1}$ (see \eqref{eq:defenergylevel} with \;$p$\; replaced by \;$p_0-1$).
\end{pf}

\section{Application of Bahri-Coron's barycenter technique}\label{eq:algtop}
In this section, we are going to use directly the results of the previous one to run a suitable scheme of the barycenter technique of Bahri-Coron\cite{bc}. To do so, we first introduce the notion of {\em neighborhood of potential critical points at infinity} of the Euler-Lagrange functional \;$\mathcal{E}_g$. Precisely, for $p\in \N^*$, $0<\varepsilon\leq \varepsilon_0$ (where $\varepsilon_0$ is given by Proposition \ref{eq:baryest}), we define \;$V(p, \varepsilon)$\; the \;$(p, \varepsilon)$-neighborhood of potential critical points at infinity of \;$\mathcal{E}_g$, namely
\begin{equation}\label{eq:vqs}
\begin{split}
V(p, \varepsilon):=\{u\in W^{1, 2}(\ov M):\;\;\exists a_1, \cdots, a_{p}\in \partial M,\;\;\alpha_1, \cdots, \alpha_{p}>0, \;\l_1, \cdots,\l_{p}\geq \frac{1}{\varepsilon},\\ ||u-\sum_{i=1}^{p}\alpha_i\varphi_{a_i, \l_i}||\leq \varepsilon,\;\;\; \;\;\frac{\alpha_i}{\alpha_j}\leq \nu_0\;\;\;\;\;\;\text{and}\;\;\;\;\;\varepsilon_{i, j}\leq \varepsilon, \;\;i\neq j=1, \cdots, p\},
\end{split}
\end{equation}
where $||\cdot||$\; denotes the standard \;$W^{1, 2}$-norm, $\varepsilon_{i, j}:=\varepsilon_{i, j}(A, \bar \l)$\; with \;$A:=(a_1, \cdots, a_{p})$, $\bar \l:=(\l_1, \cdots, \l_{p})$, \; $(\varepsilon_{i, j}(A, \bar \l))$'s are defined by \eqref{eq:varepsilonij}, $\varphi_{a_i, \l}$\; is given by \eqref{eq:varphialb} for $i=1,\cdots, p$, and \;$\nu_0$\; is given by Proposition \ref{eq:baryest}. 
\vspace{6pt}

\noindent
Concerning the sets \;$V(p, \varepsilon)$, for every \;$p\in \N^*$, we have that there exists \;$0<\varepsilon_p\leq\varepsilon_0$\; such that for every \;$0<\varepsilon\leq \varepsilon_p$, we have that
\begin{equation}\label{eq:mini}
\begin{split}
\forall u\in V(p, \varepsilon),\;\;\; \text{the minimization problem }\;\;\;\;\;\min_{B_{C_1\varepsilon}^{p}}||u-\sum_{i=1}^{p}\alpha_i\varphi_{a_i, \l_i}||
\end{split}
\end{equation}
has a unique solution, up to permutations, where \;$B^{p}_{C_1\varepsilon}$\; is defined as follows
\begin{equation}
\begin{split}
B_{C_1\varepsilon}^{p}:=\{(\bar\alpha, A, \bar \l)\in \R^{p}_+\times (\partial M)^p\times (0, +\infty)^{p}:\;\;\l_i\geq \frac{1}{\epsilon}, i=1, \cdots, p,\\\frac{\alpha_i}{\alpha_j}\leq \nu_0\;\; \text{and}\;\ \varepsilon_{i, j}\leq C_1\varepsilon, i\neq j=1, \cdots, p\}.
\end{split}
\end{equation}
and \;$C_1>1$. Furthermore, we define the selection map \;$s_{p}: V(p, \varepsilon)\longrightarrow (\partial M)^p/\sigma_p$ as follows
\begin{equation}\label{eq:select}
s_{p}(u):=A, \;\;u\in V(p, \varepsilon), \;\,\text{and} \,\;A\;\;\text{is given by}\;\,\eqref{eq:mini}.
\end{equation}
\vspace{6pt}

\noindent
Now having introduced the neighborhoods of potential critical points at infinity of the Euler-Lagrange functional \;$\mathcal{E}_g$, we are ready to present our algebraic topological argument for existence. In order to do that, we start by the following classical deformation Lemma which follows from the same arguments as for its counterparts in classical application of the algebraic topological argument of Bahri-Coron\cite{bc}(see for example Proposition 6 in \cite{bc} or Lemma 17 in \cite{BB}) and the fact that the \;$\varphi_{a, \l}$ can replace the standard bubbles in the analysis of diverging PS sequences of the Euler-Lagrange functional\; $\mathcal{E}_g$.
\begin{lem}\label{eq:classicdeform}
Assuming that \;$\mathcal{E}_g$ \;has no critical points, then for every \;$p\in \N^*$, up to taking $\varepsilon_p$\; smaller (where \;$\varepsilon_p$\; is given by \eqref{eq:mini}), we have that for every\; $0<\varepsilon\leq \varepsilon_p$, there holds \;$(W_p,\; W_{p-1})$\; retracts by deformation onto \;$(W_{p-1}\cup A_p, \;W_{p-1})$\; with \;$V(p, \;\tilde \varepsilon)\subset A_p\subset V(p, \;\varepsilon)$\; where \;$0<\tilde \varepsilon<\frac{\varepsilon}{4}$\; is a very small positive real number and depends on $\varepsilon$.
\end{lem}
\vspace{6pt}

\noindent
Using Lemma \ref{eq:baryest} and Lemma \ref{eq:classicdeform}, we are going to show that if \; $\mathcal{E}_g$\; has no critical points, then for \;$\l$\; large enough, the map \;$(f_1(\l))_*$\; is  well defined and maps \;$\partial M$\;(in top homology) in a nontrivial way in \;$(W_1,\; W_0)$. Precisely, we show:
\begin{lem}\label{eq:nontrivialf1}
Assuming that \;$\mathcal{E}_g$\; has no critical points and \:$0<\varepsilon\leq  \varepsilon_1$ (where \;$\varepsilon_1$\; is given by \eqref{eq:mini}), then up to taking \;$\varepsilon_1$\; smaller and \;$\l_1$\; larger (where \;$\l_1$\; is given by Proposition \ref{eq:baryest}), we have that for every \;$\l\geq \l_1$, there holds
$$
f_1(\l): \;(B_1(\partial M),\; B_0(\partial M))\longrightarrow (W_1, \;W_0)
$$
is well defined and satisfies
$$
(f_1(\l))_*(w_1)\neq 0\;\;\;\;\text{in}\;\;\;\;H_{n-1}(W_1, \;W_0).
$$
\end{lem}
\begin{pf}
It follows from the selection map \;$s_{1}$\; given by \eqref{eq:select}, Proposition \ref{eq:baryest} and the same arguments as in Lemma 26 in \cite{gam2}.
\end{pf}
\vspace{6pt}

\noindent
Next, like in \cite{martndia1}, using Lemma \ref{eq:transfert}, Proposition \ref{eq:baryest}, Lemma \ref{eq:classicdeform}, and the algebraic topological argument of Bahri-Coron\cite{bc}, we are going to show that if for  \;$\l$ \;large  \;$B_p(\partial M)$\; (in top homology) survives ``topologically`` the embedding into \;$(W_p, \;W_{p-1})$ via $f_p(\l)$, then for $\l$ large \;$B_{p+1}(\partial M)$\; (in top homology and as a cone with base \;$B_{p-1}(\partial M)$\; and top \;$\partial M$)  survives ''topologically`` the embedding into \;$(W_{p+1},\; W_p)$ \;via $f_{p+1}(\l)$. Precisely, we prove the following proposition:
\begin{pro}\label{eq:nontrivialrecursive}
Assuming that \;$\mathcal{E}_g$\; has no critical points and $0<\varepsilon\leq  \varepsilon_{p+1}$ (where \;$\varepsilon_{p+1}$\; is given by \eqref{eq:mini}), then up to taking \;$\varepsilon_{p+1}$\; smaller, and \;$\l_p$\; and \;$\l_{p+1}$\; larger\;(where \;$\l_p$\; and \;$\l_{p+1}$\; are given by Proposition \ref{eq:baryest}), we have that for every \;$\l\geq \max\{\l_p, \l_{p+1}\}$, there holds
$$
f_{p+1}(\l): (B_{p+1}(\partial M),\; B_{p}(\partial M))\longrightarrow (W_{p+1}, \;W_{p})
$$
and 
$$
f_p(\l): (B_p(\partial M), \;B_{p-1}(\partial M))\longrightarrow (W_p, \; W_{p-1})
$$
are well defined and satisfy
$$(f_p(\l))_*(w_p)\neq 0\;\;\;\; \text{in}\;\; \;\;H_{np-1}(W_p, \;W_{p-1})$$ implies
$$(f_{p+1}(\l))_*(w_{p+1})\neq 0\;\;\;\; \text{in} \;\;\;\;H_{n(p+1)-1}(W_{p+1}, \;W_{p}).$$
\end{pro}
\begin{pf}
First of all, we let \;$p\in \N^*$\; and \;$0<\varepsilon_{p+1}$, where \;$\varepsilon_{p+1}$\; is given by \eqref{eq:mini}. Next, recalling that we have assumed that \;$\mathcal{E}_g$\; has no critical points, and using Lemma \ref{eq:classicdeform}, then up to taking \;$\varepsilon_{p+1}$\; smaller, we infer that the following holds
\begin{equation}\label{eq:idents4}
(W_{p+1}, \;W_{p})\simeq (W_{p}\cup \mathcal{A}_{p+1}, \;W_{p}),
\end{equation}
with
\begin{equation}\label{eq:infinitynontrivial}
 V(p+1, \;\tilde \varepsilon)\subset \mathcal{A}_{p+1}\subset V(p+1,\; \varepsilon), \;\;0<4\tilde\varepsilon<\varepsilon.
\end{equation}
Now, using Lemma \ref{eq:interactestl} and Proposition \ref{eq:baryest}, we have that for every \;$\l\geq \max\{\l_p, \l_{p+1}\}$ (where $\l_p$ and $\l_{p+1}$\; are given by Proposition \ref{eq:baryest}), there holds
\begin{equation}\label{eq:mappingfp1}
f_{p+1}(\l): (B_{p+1}(\partial M), \;B_{p}(\partial M))\longrightarrow (W_{p+1}, \;W_{p}),
\end{equation}
and
\begin{equation}\label{eq:mappingfp}
f_p(\l): (B_p(\partial M), \; B_{p-1}(\partial M))\longrightarrow (W_p, \;W_{p-1}),
\end{equation}
are well defined and hence have that the first point is proven. Next, using Proposition \ref{eq:baryest}, \eqref{eq:mappingfp1}, and \eqref{eq:mappingfp}, we have that up to taking \;$\l_{p+1}$\; and \;$\l_{p}$\; larger (for example larger than \;$4\max\{\l_{p+1}(\varepsilon),\;\l_{p}(\varepsilon), \;\l_{p}(2\tilde\varepsilon), \l_p(\frac{\tilde\varepsilon}{2}), \frac{1}{\tilde \varepsilon}\}$, where $\l_p(\varepsilon)$, \;$\l_{p+1}(\varepsilon), \l_p(2\tilde\varepsilon)$,\;and \;$\l_{p}(\frac{\tilde\varepsilon}{2})$\; are given by Proposition \ref{eq:baryest} and \;$\tilde \varepsilon$ \; is given by \eqref {eq:infinitynontrivial}) the following diagram 
\begin{equation}\label{eq:diag2}
\begin{CD}
 (B_{p+1 }(\partial M), \;\mathcal{O}(B_{p}(\partial M)))  @>f_{p+1}(\lambda)>> (W_{p+1}, \; W_{p})\\
@AAA @AAA\\
(\mathcal{O}(B_{p}(\partial M)),\;B_{p-1}(\partial M))  @>f_{p}(\l)>> (W_{p}, \; W_{p-1})
\end{CD}
\end{equation}
is well defined and commutes, where
\begin{equation}\label{eq:nlbary}
\mathcal{O}(B_{p}(\partial M)):=\{\sigma=\sum_{i=1}^{p+1}\alpha_i\d_{a_i}\in B_{p+1}(\partial M):\;\;\exists \;i_0\neq j_0:\;\frac{\alpha_{i_0}}{\alpha_{j_0}}>\nu_0\;\;\;\;\text{or}\;\;\;\;\sum_{i\neq j}\varepsilon_{i, j}>\tilde\varepsilon\},
\end{equation}
with \;$\nu_0$\; given by Proposition \ref{eq:baryest}. On the other hand, we have
\begin{equation}\label{eq:isodel1}
B_{p+1}(\partial M)\setminus \mathcal{O}(B_{p}(\partial M))\simeq B_{p+1}(\partial M)\setminus B_{p}(\partial M),
\end{equation}
and
\begin{equation}\label{eq:isodel2}
\mathcal{O}(B_{p}(\partial M))\simeq  B_{p}(\partial M).
\end{equation}
Now, using \eqref{eq:purem}, Lemma \ref{eq:transfert}, and \eqref{eq:isodel1}, we derive
\begin{equation}\label{eq:tranfert2}
\begin{CD}
 H^{n-1}(B_{p+1}(\partial M)\setminus \mathcal{O}(B_{p}(\partial M)))\times H_{n(p+1)-1}(B_{p+1}(\partial M), \;B_{p}(\partial M))@>\frown>> H_{n(p+1)-n}(B_{p+1}(\partial M),\; B_{p}(\partial M))\\@>\partial>>H_{n(p+1)-n-1}(B_{p}(\partial M),\; B_{p-1}(\partial M)).
 \end{CD}
\end{equation}
Furthermore, using \eqref{eq:idents4}, we infer that
\begin{equation}\label{eq:transfert3}
\begin{CD}
 H^{n-1}(\mathcal{A}_{p+1})\times H_{n(p+1)-1}(W_{p+1},\; W_{p})@>\frown>> H_{n(p+1)-n}(
 W_{p+1},\; W_{p})\\@>\partial>>H_{n(p+1)-n-1}(W_{p},\; W_{p-1}).
 \end{CD}
\end{equation}
Moreover, passing to homologies in \eqref{eq:diag2} and using \eqref{eq:isodel2}, we derive that
\begin{equation}\label{eq:diag3a}
\begin{CD}
 H_{n(p+1)-1}(B_{p+1}(\partial M), \;B_{p}(\partial M))   &@>(f_{p+1}(\l))_*>>  H_{n(p+1)-1}(W_{p+1}, \; W_{p})
\end{CD}
\end{equation}
and
\begin{equation}\label{eq:diag3b}
\begin{CD}
 H_{np-1}(B_{p}(\partial M), \;B_{p-1}(\partial M))  & @>(f_{p}(\l))_*>>  H_{np-1}(W_{p}, \;W_{p-1})
\end{CD}
\end{equation}
are well defined and the following diagram commutes
\begin{equation}\label{eq:diag3}
\begin{CD}
 H_{n(p+1)-n}(B_{p+1}(\partial M), \;B_{p}(\partial M))   &@>(f_{p+1}(\l))_*>>  H_{n(p+1)-n}(W_{p+1},\; W_{p})\\
@V {\partial} VV &@ V {\partial} VV\\
H_{np-1}(B_{p}(\partial M), \;B_{p-1}(\partial M))  & @>(f_{p}(\l))_*>>  H_{np-1}(W_{p}, \;W_{p-1}).
\end{CD}
\end{equation}
Next, recalling that we have taken \;$\l_{p+1}$\; and \;$\l_{p}$\; larger  than \;$4\max\{\l_{p+1}(\varepsilon),\;\l_{p}(\varepsilon), \;\l_{p}(2\tilde \varepsilon), \l_p(\frac{\tilde\varepsilon}{2}), \frac{1}{\tilde \varepsilon}\}$, we derive that
\begin{equation}\label{eq:insidevm}
f_{p+1}(\l)\left(B_{p+1}(\partial M)\setminus \mathcal{O}(B_{p}(\partial M))\right)\subset V(p+1, \;\tilde \varepsilon)\subset \mathcal{A}_{p+1}\subset V(p+1, \varepsilon).
\end{equation}
Thus, using \eqref{eq:defom1}, \eqref{eq:select} and \eqref{eq:insidevm}, we infer that 
\begin{equation}\label{eq:rbc}
(f_{p+1}(\l))^*(s_{p+1}^*(O^*_{\partial M}))=O^*_{\partial M}\;\;\;\text{with}\;\; \;s^*_{p+1}(O^*_{\partial M})\neq 0\;\;\;\text{in}\;\;\; H^{n-1}(\mathcal{A}_{p+1}).
\end{equation}
On the other hand, using \eqref{eq:diag3}, we derive that
\begin{equation}\label{eq:commutation}
\partial(f_{p+1}(\l))_{*} =  (f_{p}(\l))_{*}\partial \;\;\;\text{in}\;\;\;H_{n(p+1)-n}(B_{p+1}(\partial M), \;
B_{p}(\partial M)).
\end{equation}
Now, combining \eqref{eq:purem}, Lemma \ref{eq:transfert}, \eqref{eq:isodel1}, \eqref{eq:tranfert2}, \eqref{eq:transfert3}, \eqref{eq:diag3a}, \eqref{eq:rbc}, and \eqref{eq:commutation}, we obtain
\begin{equation}\label{eq:transfim}
\begin{split}
(f_{p}(\l))_*(\o_{p})&=(f_{p}(\l))_*\left(\partial(O^{*}_{\partial M} \smallfrown \o_{p+1}) \right)\\
&=(f_{p}(\l))_*\left(\partial(((f_{p+1}(\l))^*(s_{p+1}^{*}(O^{*}_{\partial M}))) \smallfrown \o_{p+1}) \right)\\
&=\partial\left((f_{p+1}(\l))_*(((f_{p+1}(\l))^*(s_{p+1}^{*}(O^{*}_{\partial M})) \smallfrown \o_{p+1}) \right)\\
&=\partial(s_{p+1}^{*}(O^{*}_{\partial M}) \smallfrown ((f_{p+1}(\l))_*(\o_{p+1}))),
\end{split}
\end{equation}
with all the equalities holding in \;$H_{np-1}(W_{p}, \;W_{p-1})$. Hence, clearly, \eqref{eq:transfim} and the assumption 
$$(f_p(\l))_*(w_p)\neq 0\;\; \;\text{in}\; \;\;H_{np-1}(W_p, \;W_{p-1})$$ implies
$$(f_{p+1}(\l))_*(w_{p+1})\neq 0\;\;\; \text{in} \;\;\;H_{n(p+1)-1}(W_{p+1}, \;W_{p}),$$
as desired, thereby completing the proof of Proposition \ref{eq:nontrivialrecursive}.
\end{pf}
\vspace{6pt}

\noindent
Now, we are ready to present the proof of Theorem \ref{eq:existence}.\\\\
\begin{pfn}{ of Theorem \ref{eq:existence}}\\
Like in \cite{martndia1}, it follows by a contradiction argument from Corollary \ref{eq:largep}, Lemma \ref{eq:nontrivialf1} and Proposition \ref{eq:nontrivialrecursive}.
\end{pfn}

\end{document}